%&amsppt 
 
% Runs without problems under 
%    TeX 3.14159 and AmS-TeX 2.1 with amsppt.sty version 2.1d
%  The following two lines may need to be modified for non-US paper
% sizes and/or other production modes.  I chose this magnification 
% for readability

\magnification=\magstephalf
\hoffset=0.5 truein

\input amstex
\documentstyle{amsppt}

\refstyle{A}
\loadeusm
\loadbold

\define\cF{{\Cal F}}

\define\cL{{\Cal L}}\define\cP{{\Cal P}}
\define\cQ{{\Cal Q}}

\define\bbC{{\Bbb C}}\define\bbD{{\Bbb D}}
\define\bbP{{\Bbb P}}\define\bbR{{\Bbb R}}
\define\bbS{{\Bbb S}}\define\bbZ{{\Bbb Z}}

\define\eug{\operatorname{\frak g}}

\define\eup{\operatorname{\frak p}}
\define\eusl{\operatorname{\frak sl}}
\define\SO{\operatorname{SO}}
\define\SL{\operatorname{SL}}

\define\Lie{\operatorname{\eusm L}}

\define\ts{\textstyle }
\define\>{\medspace}

\define\w{{\mathchoice{\mathbin{\scriptstyle\wedge}}
                      {\mathbin{\scriptstyle\wedge}}
                      {\mathbin{\scriptscriptstyle\wedge}}
                      {\mathbin{\scriptscriptstyle\wedge}}}}

\define\smcap{{\mathchoice{\mathbin{\scriptstyle\cap}}
                          {\mathbin{\scriptstyle\cap}}
                          {\mathbin{\scriptscriptstyle\cap}}
                          {\mathbin{\scriptscriptstyle\cap}}}}

\define\lb{ [\![ }
\define\rb{ ]\!] }

\define\abold{{\bold a}} \define\bbold{{\bold b}} 
\define\gbold{{\bold g}}  \define\pbold{{\bold p}}
\define\qbold{{\bold q}}  \define\ubold{{\bold u}}
\define\vbold{{\bold v}} \define\wbold{{\bold w}} \define\xbold{{\bold x}}
\define\ybold{{\bold y}} \define\Ibold{{\bold I}} \define\Jbold{{\bold J}}
\define\Tbold{{\bold T}} \define\Wbold{{\bold W}} \define\Xbold{{\bold X}}

\define\phibold{{\boldsymbol\phi}} \define\rhobold{{\boldsymbol\rho}}
\define\zetabold{{\boldsymbol\zeta}}

\topmatter
\title 
Projectively Flat Finsler 2-spheres \\
       of Constant Curvature
\endtitle

\rightheadtext{Projectively Flat Finsler Spheres}

\author
Robert L. Bryant
\endauthor
\affil
Duke University
\endaffil
\address
Department of Mathematics,
Duke University,
PO Box 90320,
Durham, NC 27708-0320
\endaddress
\email
bryant\@math.duke.edu
\endemail

\dedicatory
\enddedicatory

\date
November 25, 1996
\enddate

\thanks
Research for this article was begun while I was visiting IMPA in Rio de
Janeiro in August 1996, while the actual writing took place while I was 
visiting the IAS in Princeton in the Fall of 1996.  Their hospitality 
as well as NSF support through the grant DMS--9505125 is hereby gratefully 
acknowledged.
\endthanks



\abstract
After recalling the structure equations of Finsler structures on surfaces, 
I define a notion of `generalized Finsler structure' as a way of 
micro-localizing the problem of describing Finsler structures subject to 
curvature conditions.  I then recall the basic notions of path geometry 
on a surface and define a notion of `generalized path geometry' analogous 
to that of `generalized Finsler structure'. 

I use these ideas to study the geometry of Finsler structures
on the $2$@-sphere that have constant Finsler-Gauss curvature~$K$ and whose 
geodesic path geometry is projectively flat, i.e., locally equivalent 
to that of straight lines in the plane. 

I show that modulo diffeomorphism there is a $2$@-parameter family of 
projectively flat Finsler structures on the sphere whose Finsler-Gauss
curvature~$K$ is identically~$1$.
\endabstract

\endtopmatter

\document

\head 0. Introduction \endhead

Hilbert's Fourth Problem was entitled ``Problem of the straight line
as the shortest distance between two points''.  It concerned, in its
most general formulation, the problem of characterizing the 
not-necessarily-symmetric distance functions~$d$ that could be defined 
on (convex) subsets~$U\subset\bbR^2$ so that the lines were geodesics, 
i.e., so that $d(x,z) \le d(x,y)+d(y,z)$ with equality if and only if 
$x$,~$y$,~and~$z$ are collinear, with~$y$ lying on the segment joining~$x$ 
to~$z$. 

Hilbert's reason for considering non-symmetric distances was
that interesting non-symmetric examples had already been 
discovered by Minkowski.  He was also aware that notions of length
of curves defined in many calculus of variations problems leads
naturally to non-symmetric distance functions.  For some interesting 
examples of physical relevance, see Carath\'eodory's book~\cite{Cara}.

The Fourth Problem can be regarded as a fundamental example of the
inverse problem in the calculus of variations.  That is, given that the
straight lines are the extremals (i.e., geodesics) of a first order 
Lagrangian for oriented curves in the plane, what can one say about 
the Lagrangian?  Can the set of such Lagrangians be characterized in
some useful way?  Can their stability properties be understood?

Beltrami had shown in 1866 that any Riemannian
metric on an open subset of the plane whose geodesics are the straight
lines must have constant Gauss curvature.  In fact, such a metric
must be locally equivalent via a projective transformation to one
of the standard metrics of constant curvature defined on the plane or
an open subset thereof, such as the Klein model of the hyperbolic plane
as the unit disk, with the geodesics being the chords of the boundary circle.
Thus, the problem was already solved for Lagrangians that represented 
Riemannian arc length.

More general Lagrangians for curves on a domain~$M$ in the plane
(or, more generally, on any surface) can be described in terms of 
 Finsler norms, where a {\it Finsler norm\/} or {\it metric\/} is a 
non-negative function~$L:TM\to\bbR$ that is positive and smooth away 
from the zero section and has the {\it homogeneity property\/} that 
$L(x,\lambda v)=\lambda\,L(x,v)$ for all~$\lambda\ge0$ and all~$v\in T_xM$ 
as well as the {\it convexity property\/} that the unit sphere
(or, as it it classically known, the {\it indicatrix\/})
$$
\Sigma_x = \{\,v\in T_xM\,\mid\, L(x,v)=1\,\}\subset T_xM
$$
at~$x\in M$ be a smooth, closed, strictly convex curve in~$T_xM$
for all~$x\in M$.  If, in addition, $L(x,-v)=L(x,v)$, or, equivalently,~
$\Sigma_x = -\Sigma_x$, then~$L$ is said to be {\it symmetric\/}.

The homogeneity property implies that for any immersed, oriented curve
~$\gamma:[a,b]\to M$ the integral
$$
\cL(\gamma) = \int_a^b L\bigl(\gamma(t),\gamma'(t)\bigr)\,dt
$$
is unchanged by oriented reparametrization and so defines a notion of
$L$@-length for oriented curves in~$M$.  The oriented curves that are 
extremal for the functional~$\cL$ are known as {\it $L$@-geodesics\/}.
When~$L$ is symmetric, the $L$@-length of a curve is independent of 
its orientation.

The convexity property implies that the Euler-Lagrange equations for the 
geo\-desics of~$L$ are everywhere non-degenerate, so that, in each oriented 
direction through each point of~$M$ there will pass a unique $L$@-geodesic 
and these $L$@-geodesics will depend smoothly on parameters, just as in the 
Riemannian case. 

The most familiar example of a Finsler norm is~$L=\sqrt g$ where~$g$ is
a Riemannian metric on~$M$.  In that case, the first
geometric invariant that can be attached to~$L$ after the geodesics 
themselves is the Gauss curvature~$K$, which is a function on~$M$.
Besides being a diffeomorphism invariant of the metric~$g$, it enters
fundamentally into the study of the second variation of the geodesics:
 For any unit speed geodesic~$\gamma:[a,b]\to M$, the Jacobi equation
for variation of geodesics near~$\gamma$ is
$$
y''(s) + K\bigl(\gamma(s)\bigr)\,y(s) = 0.
$$
 For example, $\gamma$ is locally minimizing when the operator~
$Jy=y''+K{\circ}\gamma\,y$ has zero nullity and index.

In the more general Finsler case, the second variation formula still
makes sense, but now~$K$ has to be defined, not on~$M$, but on the
unit sphere bundle~$\Sigma\subset TM$.  Then, for a unit speed geodesic~
$\gamma$ as above, the Jacobi equation turns out to be
$$
y''(s) + K\bigl(\gamma'(s)\bigr)\,y(s) = 0,
$$
and it bears the same relation to minimizing properties of~$L$@-geodesics
as the classical Jacobi equation in the Riemannian case.  This form
of the second variation was first discovered in~1907 by Underhill~\cite{Un}. 
The function~$K$ (to be described more explicitly in~\S1 below) was given 
the name `inneres Kr\"ummungsma\ss' by Finsler in his 1918 dissertation on 
the geometry of Finsler spaces, but in the present article, it will be 
called the {\it Finsler-Gauss curvature\/}.

It turned out that there were many solutions to Hilbert's Fourth Problem,
and the interested reader can profitably consult the survey article~\cite{Bu}
or the book~\cite{Po}, where a complete solution in the case of symmetric
distances is presented.  Motivated by the central place that metrics
of constant curvature play in Riemannian geometry and the direct relation
that it bears with the stability question, it is natural to pay special
attention to the Finsler norms for which~$K$ is constant.

In 1929, Funk~\cite{Fu1} showed that the Finsler examples constructed by 
Hilbert on arbitrary convex domains in the plane satisfied~$K=-1$ and
went on to classify all such solutions to Hilbert's Fourth Problem
on convex domains in the plane.  He also classified the solutions with~$K=0$. 
In each case, the local solutions depended on arbitrary functions of one 
variable.  In 1935~\cite{Fu2}, he returned to the subject and found a local 
formula for the solutions to Hilbert's Fourth Problem that satisfied~$K=+1$. 
He showed that the local solutions could be described in terms of an arbitrary 
holomorphic function of one variable subject to some inequalities. 
In 1963~\cite{Fu3}, he showed that the only global solutions on the entire 
projective plane that satisfied certain other local conditions, including 
symmetry, were the Riemannian solutions.

There remains the question of whether there are any non-Riemannian
solutions to Hilbert's Fourth Problem defined on the entire projective plane 
or the sphere and satisfying~$K=1$.  In his 1976 survey article, Busemann~
\cite{Bu, p.~139} wrote ``\dots Funk [determined] all two-dimensional 
Desarguesian metrics with constant positive curvature.  Even [these] depend on 
an essentially arbitrary function, where `essential' means that the function 
is restricted by an equality only.''   (In this context, `Desarguesian metric'
can be taken to mean `Finsler metric with linear geodesics'.)  From the context
of Busemann's remark, it seems clear that he thought that the Finsler metrics 
on the 2-sphere having constant curvature and solving Hilbert's Fourth Problem 
would depend on an `arbitrary' function.

In a 1986 paper providing a new solution to Hilbert's Fourth Problem
in the symmetric case, Szab\'o states~\cite{Sz, pp.~297--299} that he will 
discuss the non-symmetric case with constant positive curvature in the sequel, 
but, to my knowledge, this sequel has never appeared. 

In this article, I show that there are indeed non-Riemannian global solutions
to Hilbert's Fourth Problem having constant curvature~$K=1$.  Despite
Busemann's claim, I further show that, up to diffeomorphism, there is
exactly a 2@-parameter family of inequivalent solutions and that the only
symmetric solution is the Riemannian one.  An explicit formula for these
 Finsler metrics is given in Theorem~10.  The unit sphere or indicatrix 
in these examples turns out to be a smooth algebraic plane curve of degree~4 
(or, at some special points or in the Riemannian case,~2) and the only 
symmetric Finsler structure among them is the Riemannian one.

This article is organized as follows:  In \S1, I recall the basics
of the geometry of Finsler structures and introduce a generalization
that will be needed in the rest of the article.  This generalization
allows for multi-valued Finsler norms (see Example~1) and other apparent
pathologies, such as Finsler structures on orbifolds, but it has the virtue
that the differential equations for prescribed curvature in Finsler
geometry are local on the generalized objects, while they are not on
the classical objects.  This generalization is actually of independent
interest and may have applications to control theory.

In \S2, I go over the basics of path geometry, i.e., the geometry of a 
surface endowed with a 2@-parameter family of paths, such as the 
geodesic path geometry of a Finsler metric.  I introduce a notion of
generalized path geometry that corresponds to the earlier notion of 
generalized Finsler structure and use this to develop a notion of and a
test for projective flatness for a generalized Finsler structure.  I
recall Cartan's projective connection for a path geometry and note that 
his construction extends without change to the case of a generalized
path geometry.  His computation of projective curvatures is used to 
recover Berwald's characterization of the Finsler structures whose
geodesics can be mapped to straight lines in the plane by local 
diffeomorphism (the projectively flat case) and to extend this concept 
to the generalized case.

In \S3, I classify the projectively flat Finsler structures with~$K=1$
that can be defined on compact simply-connected 3@-manifolds.  I show that
up to diffeomorphism, there is a 2@-parameter family of these and and that
they each arise as the double cover of the unit sphere bundle of a
classical Finsler structure on the 2-sphere.  The main tool in the finiteness
theorem is a vanishing theorem for a certain holomorphic cubic differential
on a Riemann surface of genus zero.  The existence theorem uses the
vanishing theorem to provide enough extra local equations to identify
the space of solutions with the leaves of an integrable distribution
of dimension~3 defined on a certain manifold~$X$ of dimension~13. 

The constructions of~\S3 are invariant under the projective group acting
on~$\bbR\bbP^2$ and they suggest that a projectively invariant interpretation 
of Funk's local characterization of these Finsler structures in the classical 
case might lead to a global classification theorem.  In~\S4, I do just this.
I recast Funk's results of~\cite{Fu1-3} as statements about sections of
a certain bundle over the double cover of~$\bbR\bbP^2$ and the geometry
of certain holomorphic curves in~$\bbC\bbP^2$.  Global considerations from
algebraic geometry make it possible to then classify the holomorphic curves
that can arise in any global solution as the conics without real points.
Then, using a classification of the moduli of conics without real points
under the action of the real group~$\SL(3,\bbR)$, I derive the explicit
formulae of Theorem~10.  I then close the paper with a few remarks about
the geometry of the examples.

Throughout the article, smoothness is assumed.  Weaker differentiability 
assumptions would have sufficed, but would have been a distraction from
the main arguments.  It appears likely that the classification results
would hold without change as long as the Finsler structures were 
differentiable enough to define the curvature, but I have not pursued
this question.

\head 1. Finsler Surfaces and Generalized Finsler Structures \endhead

This section contains a brief review of basic Finsler surface theory and
Cartan's construction of the canonical coframing associated to a Finsler 
structure on a surface.  Along with this review, I will also introduce and 
discuss a notion of generalized Finsler structure that will be used
in the rest of the article. 

The material on classical Finsler surfaces has been treated in 
many places, cf.~\cite{BaChSh}, \cite{Ca2}, \cite{Ch1}, \cite{Ru}, \cite{Ma2}, 
and \cite{GaWi}, to name just a few from different points of view and 
different eras.  It is included here to establish notation 
and nomenclature. 

Let $M$ be a connected, smooth, oriented surface (i.e., a $2$@-manifold). 
A {\it Finsler structure\/} on~$M$ 
is a smooth hypersurface~$\Sigma^3\subset TM$ for which the basepoint 
projection~$\pi:\Sigma\to M$ is a surjective submersion and 
having the property that for each~$x\in M$, the $\pi$-fiber~$\Sigma_x
=\pi^{-1}(x)=\Sigma\,\smcap\, T_xM$ is a closed, strictly convex curve 
enclosing the origin~$0_x\in T_xM$.  If, in addition, the curve~$\Sigma_x$ 
is symmetric about~$0_x$ for each~$x\in M$, then $\Sigma$ will be said to be 
{\it symmetric\/}.

A differentiable curve~$\gamma:[a,b]\to M$ will be said to be
a~{\it $\Sigma$@-curve\/} if, for every~$s$ in the interval~$[a,b]$, 
the velocity vector~$\gamma'(s)$ lies in~$\Sigma$.  The map~
$\gamma':[a,b]\to\Sigma$ is known as the {\it tangential lift\/}
of~$\gamma$.  For every immersed curve~$\gamma:[a,b]\to M$, there is
a unique orientation preserving diffeomorphism~$h:[0,L]\to[a,b]$ so
that~$\gamma\circ h$ is a $\Sigma$@-curve.  The number~$L>0$ is the
{\it $\Sigma$@-length\/} of~$\gamma$.  Note that if $\Sigma$ is 
not symmetric then reversing the orientation of an immersed curve
may change its length.

\subhead 1.1.  The canonical coframing \endsubhead
In~\cite{Ca2}, \'E.~Cartan constructed a canonical coframing on any 
 Finsler structure~$\Sigma\subset TM$ on an oriented surface~$M$.  In 
that paper, he essentially proves the following result.

\proclaim{Theorem 1 (Cartan)}
Let $\Sigma\subset TM$ be a Finsler structure on an oriented surface~$M$
and let~$\pi:\Sigma\to M$ denote the basepoint projection.  Then there 
exists a unique coframing~$\omega=\bigl(\omega_1,\omega_2,\omega_3\bigr)$ 
of~$\Sigma$ with the properties:
\roster
\item $\omega_1\w\omega_2$ is a positive multiple of any 
      $\pi$@-pullback of a positive $2$@-form on~$M$;
\item The tangential lift of any $\Sigma$-curve~$\gamma$ 
      satisfies~$(\gamma')^*\omega_2=0$ and $(\gamma')^*\omega_1 = ds$; 
\item $\omega_2\w d\omega_1=0$;
\item $\omega_1\w d\omega_1 = \omega_2\w d\omega_2$;
\item $d\omega_1 = -\omega_2\w\omega_3$ and $\omega_3\w d\omega_2 = 0$.
\endroster
Moreover, there exist unique functions~$I$, $J$, and $K$ on~$\Sigma$
so that
$$
\aligned
d\omega_1 &= {} -\omega_2\w\omega_3\,,\\
d\omega_2 &= {} -\omega_3\w\bigl(\omega_1 - I\,\omega_2\bigr)\,, \\
d\omega_3 &= {} -\bigl(K\,\omega_1 - J\,\omega_3\bigr)\w\omega_2\,.\\
\endaligned
\tag1
$$
\endproclaim

Let~$\Xbold=\bigl(\Xbold_1,\Xbold_2,\Xbold_3\bigr)$ be the vector field 
framing of~$\Sigma$ that is dual to the coframing~$\omega=
\bigl(\omega_1,\omega_2,\omega_3\bigr)$.  In terms of~$\Xbold$, equations~(1) 
can be expressed in the form
$$
\aligned
\bigl[\,\Xbold_2,\,\Xbold_3\,\bigr] &= \Xbold_1+I\,\Xbold_2+J\,\Xbold_3,\\
\bigl[\,\Xbold_3,\,\Xbold_1\,\bigr] &= \Xbold_2,\\
\bigl[\,\Xbold_1,\,\Xbold_2\,\bigr] &= K\,\Xbold_3.\\
\endaligned
\tag 2
$$

Note that the fibers of the basepoint projection~$\pi$ are simply the 
integral curves of the vector field~$\Xbold_3$.

\remark{Remarks}
The notation is Cartan's. However, the reader should be aware that 
other authors have modified it somewhat. For example, H. Rund denotes the 
invariant~$I$ by the letter~$J$. In a recent paper~\cite{Br}, I used the 
letters~$S$ and $C$ to denote the quantities $-I$ and $-J$, a decision that 
I now regret.

Strictly speaking, there is a slight difference between Theorem~1
and the result of Cartan.  Since Cartan does not fix an orientation on
the surface~$M$, the first property in the list above has no meaning for him.
As a result, Cartan's coframing is only well-defined up to an ambiguity
of the form~$\omega = \bigl(\omega_1,\pm\omega_2,\pm\omega_3\bigr)$. 
The introduction of an orientation on~$M$ fixes this ambiguity.
 
In standard treatments of the calculus of variations, 
the $1$@-form~$\omega_1$ is known as {\it Hilbert's invariant integral\/}.

A $\Sigma$-curve~$\gamma$ is a {\it $\Sigma$@-geodesic\/}, i.e., a critical 
point of the natural $\Sigma$@-length functional on curves in~$M$, if and only 
if its tangential lift satisfies~$(\gamma')^*\omega_3=0$. Thus, the 
$\Sigma$@-geodesics are the projections to~$M$ of the integral curves 
of~$\Xbold_1$.  For this reason, the flow of the vector field~$\Xbold_1$
is known as the {\it geodesic flow} of~$\Sigma$ and $\Sigma$ is
said to be {\it geodesically complete\/} if ~$\Xbold_1$ is complete, i.e., its 
flow from any initial point is defined for all time (positive and negative).

The function~$I$ vanishes if and only if~$\Sigma$ is the unit circle bundle
of a Riemannian metric~$g$ on~$M$.  In this case, the function~$K$ is
the $\pi$@-pullback to~$\Sigma$ of the Gauss curvature function of~$g$.
\endremark

\subhead 1.2.  Generalized Finsler structures \endsubhead
Taking Cartan's construction as a starting point leads to a natural
widening of the notion of a Finsler structure. 

\definition{Definition 1}
A {\it generalized Finsler structure} on a $3$@-manifold~$\Sigma$ is a 
coframing~$\omega=\bigl(\omega_1,\omega_2,\omega_3\bigr)$ that satisfies the 
equations~$(1)$ for some (necessarily unique) functions~$I$, $J$, and $K$ 
on~$\Sigma$.
\enddefinition

As in the classical case, for a generalized Finsler structure~$\omega=
\bigl(\omega_1,\omega_2,\omega_3\bigr)$, the dual framing of vector fields 
will be denoted~$\Xbold=\bigl(\Xbold_1,\Xbold_2,\Xbold_3\bigr)$. 

A generalized Finsler structure will be said to be {\it amenable\/} if the 
leaf space of the foliation defined by the integral curves of~$\Xbold_3$ can 
be given the structure of a smooth surface~$M$ in such a way that the natural 
projection~$\pi:\Sigma\to M$ is a smooth submersion.

Every generalized Finsler structure~$\bigl(\Sigma,\omega\bigr)$ is locally 
amenable in the sense that every point of~$\Sigma$ has a neighborhood to 
which the generalized Finsler structure restricts to be amenable.  In fact, 
the next proposition shows that the difference between the concepts 
`Finsler structure' and `generalized Finsler structure' is global in 
nature; every generalized Finsler structure is locally diffeomorphic 
to a Finsler structure.  The proof is straightforward, so I omit it.

\proclaim{Proposition 1}
Let~$\Sigma$ be a $3$-manifold endowed with an amenable generalized Finsler
structure~$\omega$.  Denote the projection onto the space of integral
curves of~$\Xbold_3$ by~$\pi:\Sigma\to M$ and define a smooth map~
$\nu:\Sigma\to TM$ by the rule~
$\nu(\ubold)=\pi'(\ubold)\bigl(\Xbold_1(\ubold)\bigr)$ for all~
$\ubold\in\Sigma$. 
Then $\nu$ immerses each $\pi$-fiber~$\Sigma_x = \pi^{-1}(x)$ as a curve 
in~$T_xM$ that is strictly convex towards~$0_x$.  Moreover, there is
an orientation of~$M$ so that the $\nu$@-pullback of the canonical 
coframing induced on the $\nu$@-image of~$\Sigma$ co\"\i ncides with 
the given generalized Finsler structure. \qed
\endproclaim

The reader may well wonder why anyone would bother with generalized Finsler
structures since they are locally the same as Finsler structures.  The
reason is that in the study of Finsler structures defined by geometric
conditions, such as conditions on the invariants~$I$, $J$, and $K$, one
is frequently led to solve differential equations in the larger class of
generalized Finsler structures since it is this class that is locally defined.
Then, as a separate step, one can determine the necessary conditions
on a generalized Finsler structure that it actually be a Finsler structure. 
Thus, generalized Finsler structures provide a natural intermediate
stage where problems can be localized and solved without the complication
of global issues.

The next proposition gives a simple necessary and sufficient test 
for a generalized Finsler structure to be a Finsler structure. The
proof is straightforward.

\proclaim{Proposition 2}
A generalized Finsler structure~$\omega$ on a $3$-manifold~$\Sigma$ is a 
 Finsler structure if and only if the integral curves of~ 
$\Xbold_3$ are all closed, it is amenable, and the canonical 
immersion~$\nu:\Sigma\to TM$ is one-to-one. \qed
\endproclaim

Note that just having the integral curves of~$\Xbold_3$ be closed does not 
make a generalized Finsler structure amenable since one could have a discrete 
subset of the leaves around which the foliation is not locally a product, 
the so-called `ramified' orbits.  In this case, the leaf space will have the 
structure of a $2$@-dimensional orbifold near the ramified points. 

When a generalized Finsler structure is amenable in addition to having
all the integral curves of~$\Xbold_3$ closed, the $\nu$-image of each 
$\pi$@-fiber~$\Sigma_x$ will in general be a closed, strictly convex curve 
in~$T_xM$ which winds around the origin $\mu$ times for some positive integer~
$\mu$ which may well be greater than~$1$.  The number~$\mu$, known as the 
{\it multiplicity\/} of the generalized Finsler structure, is equal to~$1$ 
if and only if each~$\Sigma_x$ is embedded via~$\nu$.

\example{Example 1}
Let $S^2$ be given its standard Riemannian metric~$g_0$ and let~$\Sigma_0
\subset TS^2$ be its unit tangent bundle, with basepoint projection~
$\pi_0:TS^2\to S^2$.  Of course, $\Sigma_0$ defines a Finsler 
structure on~$S^2$ in the classical sense. 

Now, it is well known that $\Sigma_0$ is diffeomorphic to~$\SO(3)$ and 
hence has~$\bbZ_2$ as its fundamental group.  Let~$\tau:\Sigma\to\Sigma_0$ 
be its universal cover, so that~$\Sigma$ is diffeomorphic to~$S^3$.  Now let~
$\tilde\tau:\Sigma\to TS^2$ be an arbitrary small perturbation of~$\tau$
and define~$\pi=\pi_0\circ\tilde\tau$.  Provided this perturbation is 
sufficiently small in the $C^2$@-topology, the $\tilde\tau$@-image of each 
fiber~$\Sigma_x = \pi^{-1}(x)\subset\Sigma$ will be a closed, immersed curve 
in~$T_xS^2$ which is strictly convex towards~$0_x$. For the generic such
perturbation, these image curves will have winding number~$2$ about~$0_x$ 
but will not double cover an embedded curve of winding number~$1$ about~$0_x$. 

Now, Cartan's construction of his canonical coframing on a hypersurface
in~$TM$ is local and depends only on the assumptions that the hypersurface
submerses onto~$M$ with fibers that are locally strictly convex towards
the origin in each tangent space.  It follows that the immersion~$\tilde\tau$
induces a generalized Finsler structure~$\tilde\omega$ on~$\Sigma$. 
This generalized Finsler structure is amenable, with~$\pi:\Sigma\to 
S^2$ being the leaf projection and $\tilde\tau$ being the canonical 
immersion~$\nu$.  In this case, the multiplicity~$\mu$ is equal to~$2$.
\endexample

\subhead 1.3.  The Bianchi identities \endsubhead
The structure equations of a generalized Finsler 
structure~$\omega$ on a $3$@-manifold~$\Sigma$ are
$$
\aligned
d\omega_1 &= {} - \omega_2\w\omega_3\,,\\
d\omega_2 &= {} - \omega_3\w\bigl(\omega_1 - I\,\omega_2\bigr)\,,\\
d\omega_3 &= {} - \bigl(K\,\omega_1 - J\,\omega_3\bigr)\w\omega_2\,,\\
\endaligned 
$$
where~$I$, $J$, and $K$ are smooth functions on~$\Sigma$. Still following
Cartan, I will now derive the Bianchi identities, of the structure, i.e., the
relations among the derivatives of the invariants~$I$, $J$, and $K$. 

The exterior derivative of the first of these three equations is an
identity when the other two are taken into account.
The exterior derivative of the second equation simplifies to
$$
0 = d\bigl(d\omega_2\bigr) = -\bigl(dI-J\,\omega_1\bigr)\w\omega_2\w\omega_3,
$$
while the exterior derivative of the third equation simplifies to
$$
0 = d\bigl(d\omega_3\bigr) 
=\bigl(-dK\w\omega_1+(dJ+KI\,\omega_1)\w\omega_3\bigr)\w\omega_2.
$$
These two equations imply that there exist functions~$I_2$,~$I_3$, $J_2$,
~$J_3$,~$K_1$, $K_2$, and $K_3$ for which
$$
\alignedat4
dI ={}&&        J\,\omega_1 &&{}+ I_2\,\omega_2 &&{}+ I_3\,\omega_3\,,&{}\\
dJ ={}&&-(K_3+KI)\,\omega_1 &&{}+ J_2\,\omega_2 &&{}+ J_3\,\omega_3\,,&{}\\
dK ={}&&      K_1\,\omega_1 &&{}+ K_2\,\omega_2 &&{}+ K_3\,\omega_3\,.&{}\\
\endaligned
\tag 3
$$

The formulae~$(3)$ constitute the Bianchi identities of the structure.
Alternatively, they can be expressed as follows:  For any differentiable 
function~$F$ on~$\Sigma$, define the functions~$F_1$, $F_2$, and $F_3$ by
the equation~$dF = F_1\,\omega_1+F_2\,\omega_2+F_3\,\omega_3$.  Then the
Bianchi identities can also be expressed in the form
$$
I_1 - J = 0
\qquad\qquad\text{and}\qquad\qquad
J_1 + K_3 + KI = 0.
$$

\subhead 1.4.  Geodesics and the second variation \endsubhead
Let~$\omega=\bigl(\omega_1,\omega_2,\omega_3\bigr)$ be a generalized Finsler 
structure on a $3$@-manifold~$\Sigma$.  The {\it geodesics\/} of the structure 
are the integral curves of~$\Xbold_1$.  They define a foliation of~$\Sigma$ 
called the {\it geodesic foliation\/}. The (local) flow of~$\Xbold_1$ is the 
{\it geodesic flow\/} and $\omega$ is said to be
{\it geodesically complete\/} if~$\Xbold_1$ is complete, i.e., 
its flow exists for all time.   The structure~$\omega$ is {\it geodesically 
amenable\/} if the  leaf space~$\Lambda$ of the geodesic foliation can be 
given the structure of a smooth surface in such a way that the natural 
projection~$\ell:\Sigma\to\Lambda$ is a smooth submersion. Of course,
every generalized Finsler structure is locally geodesically amenable.

When $\Sigma\subset TM$ is a Finsler structure on a surface~$M$, a 
$\Sigma$@-geodesic is a $\Sigma$@-curve~$\gamma:D\to M$ (where $D\subset\bbR$
is some interval, which may be bounded or unbounded and closed or open) 
for which the lifted curve~$\gamma':D\to\Sigma$ is a geodesic of the 
generalized Finsler structure on~$\Sigma$. If $\gamma:D\to M$ is a 
$\Sigma$@-geodesic and the interval~$D$ contains~$0$ then 
$\gamma'(s)=\exp_{s\Xbold_1}(\ubold)$ where $\ubold=\gamma'(0)$. 
Thus, $\gamma(s) = \pi\bigl(\exp_{s\Xbold_1}(\ubold)\bigr)$.

In the Riemannian geometry of surfaces, the Gaussian curvature plays
an important role in the formula for the second variation of arc length.
In the more general case of a Finsler structure, the function~$K$,
known as the Finsler-Gauss curvature, plays
the same role.   More precisely, a compact $\Sigma$@-geodesic~
$\gamma:[a,b]\to M$ will be a local minimum of the $\Sigma$@-length functional
if the quadratic form
$$
Q_\gamma(f) = \int_a^b \bigl((f')^2 - (K{\circ}\gamma')\,f^2\bigr)\,dt
$$
has zero index and nullity on the space of smooth functions~$f$ on~$[a,b]$ 
which vanish at the endpoints.  This follows by an elementary calculation 
which I omit.

In particular, if $K$ is non-positive then every geodesic segment
is locally minimizing.  On the other hand, if $K\ge a^2$ for some positive
constant~$a$, then, just as in the Riemannian case, no geodesic segment
of length greater than $\pi/a$ can be locally minimizing.

\remark{Remark}
When a Finsler structure~$\Sigma\subset TM$ is geodesically amenable,
Crofton's Formula holds in the following sense:  If $\ell:\Sigma\to\Lambda$
is the submersion onto the leaf space of the geodesic flow, then
there exists a unique $2$@-form~$\mu$ on~$\Lambda$ so that ~$\ell^*\mu=
\omega_2\w\omega_3$.  For any immersed curve~$\gamma:[0,1]\to M$, let
$\gamma^-:[0,1]\to M$ be the same curve traversed in the opposite
orientation.  For any oriented geodesic~$\lambda\in\Lambda$ that meets~
$\gamma$ transversely, let $\nu^+_\gamma(\lambda)$ be the number of 
positively oriented intersections of~$\lambda$ with~$\gamma$.  Define
~$\nu^-_\gamma$ similarly to be the number of negatively oriented 
intersections of~$\lambda$ with~$\gamma$.  Then, as Berwald observed,
Crofton's formula remains valid in this case in the form
$$
\cL(\gamma)+\cL(\gamma^-) = \int_\Lambda \nu^+_\gamma\,\mu
= \int_\Lambda \nu^-_\gamma\,\mu\,.
$$
When~$\Sigma$ is symmetric, $\cL(\gamma)=\cL(\gamma^-)$ and the formula
simplifies to
$$
\cL(\gamma) = {\frac{1}{4}}\,\int_\Lambda (\nu^+_\gamma+\nu^-_\gamma)\,\mu\,.
$$
This shows that, in the geodesically amenable case, the notion of
$L$@-length in~$M$ can be recovered from the knowledge of the space~$\Lambda$ 
of $L$@-geodesics together with a measure~$\mu$ on~$\Lambda$. 

This observation forms the basis of Pogorelov's solution of 
Hilbert's Fourth Problem in the symmetric case, as well as Szab\'o's more 
recent treatment.  This also points out the importance of studying the path 
geometry defined by the $L$@-geodesics.  It is to this study that the
next section is devoted.
\endremark

\head 2. Path Geometries and Projective Structures \endhead

A Finsler structure on a surface~$M$ defines a 2-parameter family of 
oriented paths on~$M$, one in every oriented direction through every point.
This is a special case of the notion of {\it path geometry\/}.  In this
section, the basics of path geometries on surfaces are recalled and a
notion of {\it generalized path geometry\/} is introduced
so as to correspond with generalized Finsler structures.  Cartan's solution
of the equivalence problem for (generalized) path geometries is then reviewed
with the purpose of developing an effective test for when a path geometry
is locally equivalent to the `flat' example of lines in the plane.  Finally,
this is applied to derive the classical conditions (due to Berwald) for
a Finsler structure to have its geodesic path geometry be equivalent to that
of the lines in the plane.

\subhead 2.1.  Path geometries and generalized path geometries \endsubhead
Roughly speaking, a {\it path geometry\/} on a surface~$M$ is a $2$@-parameter
family~$\Lambda$ of curves on~$M$ with the property that for every point
$x\in M$ and every line~$L\subset T_xM$ there is a unique curve~$\xi\in\Lambda$
with the property that~$\xi$ passes through~$x$ and has~$L$ as its tangent
line at~$x$.  The fundamental example to keep in mind is the family of lines
in the Euclidean plane.  The actual definition to be given below will refine 
this intuitive picture by incorporating appropriate notions of smoothness
and independence.

\subsubhead 2.1.1.  Path geometry on surfaces \endsubsubhead
 For any surface~$M$, let $TM$, as usual, denote the tangent bundle, and
let~$\pi:\bbP(TM)\to M$ denote the projectivized tangent bundle, whose fiber 
over~$x\in M$ is the projectivization of~$T_xM$, i.e., the space of lines
through~$0_x\in T_xM$.  Given any smooth, immersed curve~$\gamma:(a,b)\to M$,
there is a canonical lift~$\gamma_1:(a,b)\to\bbP(TM)$ defined by the 
rule~$\gamma_1(t) = T_{\gamma(t)}\gamma\bigl((a,b)\bigr)$.

It is easy to characterize the canonical lifts of immersed curves in~$M$ as
curves in~$\bbP(TM)$ in terms of local geometric structures on~$\bbP(TM)$. 
In fact, as is well-known, $\bbP(TM)$ carries the structure of a 
$3$@-dimensional contact manifold.  This contact structure is defined as 
follows:  Since $\pi:\bbP(TM)\to M$ is a submersion, for each~$L\in\bbP(TM)$, 
the linear map $\pi'(L):T_L\bbP(TM)\to T_{\pi(L)}M$ is surjective.  Hence, the 
inverse image, $E(L) = \pi'(L)^{-1}(L)\subset T_L\bbP(TM)$ is a canonically 
defined $2$@-plane in~$T_L\bbP(TM)$.  The $2$@-plane field~$E$ defines a 
contact structure on~$\bbP(TM)$. 

 From the definition of the canonical lift~
$\gamma_1$ given above, it follows that $\gamma_1$ is tangent to~$E$ at all 
points, i.e., $\gamma_1$ is a contact curve.  Moreover, since the projection~
$\gamma=\pi\circ\gamma_1$ is an immersion, it follows that $\gamma_1$ is also 
transverse to the fibers of~$\pi$.  Conversely, any contact curve~$\phi:(a,b)
\to\bbP(TM)$ which is transverse to the fibers of~$\pi$ is of the form~$\phi
=\bigl(\pi\circ\phi\bigr)_1$, and so is, in particular, a canonical lift.

\definition{Definition 2}
A {\it path geometry\/} on a surface~$M$ is a foliation~$\cP$ of~$\bbP(TM)$ by
contact curves, each of which is transverse to the fibers of~$\pi:\bbP(TM)
\to M$.  A {\it local path geometry\/} on~$M$ is a foliation~$\cP$ of an open 
subset~$U\subset\bbP(TM)$ by contact curves, each of which is transverse 
to the fibers of~$\pi$.
\enddefinition

Note that, since a path geometry~$\cP$ is a codimension~2 foliation, it makes
sense to say that it defines a $2$@-parameter family of curves on the base 
space~$M$, though this is somewhat imprecise.  It could happen that the space 
of leaves of~$\cP$ is non-Hausdorff. 

In the case that there is a surface~$\Lambda$
and a submersion~$\lambda:\bbP(TM)\to\Lambda$ whose fibers are the leaves
of~$\cP$, the path geometry will be said to be {\it amenable\/}.  This is the
case, for example, for the path geometry that consists of the lines in
the standard Euclidean plane.  By contrast, the path geometry defined by the 
geodesics on a compact surface of constant negative curvature is very far
from being amenable, as the geodesic flow is ergodic in this case.

\subsubhead 2.1.2.  Generalized path geometries \endsubsubhead
In this article, it will be necessary to work with
path geometries that are only locally defined.  In fact, it is useful to
generalize the notion of a path geometry as follows:

\definition{Definition 3}
A {\it generalized path geometry\/} on a $3$@-manifold~$\Sigma$ is a
pair of transverse codimension~2 foliations~$\bigl(\cP,\cQ\bigr)$ with the 
property that the (unique) $2$@-plane field~$E$ that is tangent to both 
foliations defines a contact structure on~$\Sigma$.
\enddefinition

A path geometry~$\cP$ on~$\Sigma=\bbP(TM)$ in the classical sense is a special 
case of a generalized path geometry where the second foliation~$\cQ$ is taken
to be the fibers of the basepoint projection~$\pi:\bbP(TM)\to M$. 

It is convenient to introduce a notion of {\it orientability\/}. 
A generalized path geometry~
$\bigl(\Sigma,\cP,\cQ\bigr)$ will be said to be {\it oriented\/} if a
continuous choice of orientation of the leaves of each of~$\cP$ and $\cQ$
has been made.  A choice of orientation is then equivalent to the choice of
two non-vanishing vector fields on~$\Sigma$ defined up to positive
multiples, say, $X_1$~and~$X_3$, the first being tangent to the leaves of~$\cP$
and the second being tangent to the leaves of~$\cQ$. 

I skipped the index~`2' because I want to set $X_2 = [X_3,X_1]$ so as to 
correspond more closely to the notation established in the previous section. 
The hypothesis that the $2$@-plane field~$E$ spanned by~$X_1$ and $X_3$ be a 
contact plane field implies that the vector field~$X_2$ is linearly independent
from~$X_1$ and $X_3$. 

Note that if I replace~$X_1$ by $X_1^*=\lambda_1\,X_1$ and $X_3$ 
by~$X_3^*=\lambda_3\,X_3$, where the functions~$\lambda_1$ and $\lambda_3$
are both positive, then $X_2$ is replaced by 
$$
X_2^* = \lambda_1\lambda_3\,X_2 +d\lambda_1(X_3)\,X_1-d\lambda_1(X_1)\,X_3\,.
$$
This remark will be useful in the next subsection and farther along.

 Finally, a $1$@-form~$\alpha$ on~$\Sigma$ will be said to be 
{\it $\cP$@-positive\/} if it pulls back to each leaf of~$\cP$ to be
positive with respect to the specified orientation.  The notion of
{\it $\cQ$@-positivity\/} is defined similarly.

\example{Example 2}
Let~$\bbR^2$ be endowed with its usual Euclidean metric and let
$\Sigma\subset \bbR^2\times\bbR^2$ be the set of pairs of points~$(\xbold,
\ybold)$ in the plane satisfying~$|\xbold-\ybold|=1$.  Let~
$\pi:\Sigma\to\bbR^2$ be the projection onto the second factor and let
$\lambda:\Sigma\to\bbR^2$ be the projection onto the first factor.  Let
$\cP$ be the foliation defined by the fibers of~$\lambda$ and let ~$\cQ$
be the foliation defined by the fibers of~$\pi$.  Then it can be verified
that~$\bigl(\Sigma,\cP,\cQ\bigr)$ is a generalized path geometry.  Moreover,
it can be oriented so that the clockwise orientation of the unit circles~
$\pi^{-1}(\ybold)$ and $\lambda^{-1}(\xbold)$ are each positive.

Note, by the way, that $\Sigma$ does not induce a path geometry on~$\bbR^2$
in the usual sense because through each point in the plane there pass 
{\it two\/} unit circles having the same (unoriented) tangent line.
\endexample 

\subsubhead 2.1.3.  Local realization \endsubsubhead
In spite of Example~2, every generalized path geometry~
$\bigl(\cP,\cQ\bigr)$ on a $3$@-manifold~$\Sigma$ is locally 
identifiable with a local path geometry on a surface. 

To show this, a (local) candidate surface~$M$ must be found. 
This can be done as follows:  Let $\ubold\in\Sigma$ be
chosen and let~$U\subset\Sigma$ be an open neighborhood of~$\ubold$ on which 
the foliation~$\cQ$ is amenable, i.e., so that there exists a smooth 
surface~$M$ and a smooth surjective submersion~$\pi:U\to M$ so that the 
fibers of~$\pi$ are the leaves of~$\cQ$ restricted to~$U$. (Note that $M$ 
and~$\pi$ are uniquely determined by~$U$ up to the 
natural notion of equivalence up to diffeomorphism.)

A canonical smooth map~$\nu:\Sigma\to\bbP(TM)$ can now be defined as follows: 
 For each $\vbold\in U$, let $\nu(\vbold)=\pi'(\vbold)\bigl(T_\vbold\cP\bigr)$. 
Since the foliations~$\cQ$ and $\cP$ are transverse and since the fibers
of~$\pi$ are the leaves of~$\cQ$, it follows that 
$\pi'(\vbold)\bigl(T_\vbold\cP\bigr)$ is a one-dimensional subspace of
$T_{\pi(\vbold)}M$, i.e., an element of~$\bbP\bigl(T_{\pi(\vbold)}M\bigr)$.
Thus, $\nu$ is well-defined.  The hypothesis that the plane field~$E$ be
a contact field implies that $\nu$ is an immersion and hence, for dimension
reasons, a local diffeomorphism.  Shrinking~$U$ if necessary, it can be
assumed that~$\nu$ embeds~$U$ as an open subset of~$\bbP(TM)$.

It is now not difficult to verify that~$\nu$ is a contact diffeomorphism, i.e., 
it identifies the given contact structure on~$\Sigma$ with the one got by 
pulling back the canonical contact structure on~$\bbP(TM)$ via~$\nu$. 

In this way, the generalized path geometry on~$U$ is canonically 
identified with a local path geometry on~$M$.

If the path geometry is oriented, there is an induced orientation
on the surface~$M$, defined as follows: Let $X_1$, $X_2$ and $X_3$ be chosen 
to correspond to the given orientation as in the last section.  The vector 
fields $X_1$ and $X_2$ then span a $2$@-plane field $H$ complimentary to the 
tangents to the fibers of~$\pi$.  Since the fibers of~$\pi$ are, by definition,
connected, it follows that there is a unique orientation on the tangent spaces
to~$M$ so that, for all~$\vbold\in U$, the isomorphism~$\pi'(\vbold):H_\vbold
\to T_{\pi(\vbold)}M$ carries~$\bigl(X_1(\vbold),X_2(\vbold)\bigr)$ to an 
oriented basis of~$T_{\pi(\vbold)}M$.  Because of the way~$X_2$ was defined
and because of the formula from \S2.1.2 for how it changes when 
~$X_1$ and $X_3$ are replaced by positive multiples, it follows that this
orientation of~$M$ depends only on the orientation of the generalized path 
geometry, not on the specific choices of~$X_1$ and $X_3$.

\subsubhead 2.1.4.  Equivalence and duality \endsubsubhead
Two generalized path geometries~$\bigl(\Sigma_1,\cP_1,\cQ_1\bigr)$ and
$\bigl(\Sigma_2,\cP_2,\cQ_2\bigr)$ are said to be {\it equivalent\/}
if there exists a diffeomorphism~$\phi:\Sigma_1\to\Sigma_2$ that
satisfies~$\phi\bigl(\cP_1)=\cP_2$ and $\phi\bigl(\cQ_1)=\cQ_2$. 

The group of self-equivalences of a given generalized path geometry~
$\bigl(\Sigma,\cP,\cQ\bigr)$ is known as its {\it automorphism group\/} or 
{\it symmetry group\/}.  By a theorem of Cartan~\cite{Ca1}, this symmetry
group has the structure of a Lie group of dimension at most~$8$.  In fact,
if $\Sigma$ is connected and $\eug$ is the set of vector fields~$X$ 
on~$\Sigma$ having the property that the (locally defined) flow of~$X$ is a 
$1$@-parameter group of (local) symmetries of~$\bigl(\Sigma,\cP,\cQ\bigr)$, 
then~$\eug$ is a Lie algebra of dimension at most~$8$, with equality
if and only if the generalized path geometry is locally equivalent to the
one defined by the path geometry of straight lines in the plane.  In this
most symmetric case, the algebra~$\eug$ is isomorphic to~$\eusl(3,\bbR)$, 
the Lie algebra of infinitesimal projective transformations of the plane.
(For further discussion of this point, see the next section.)

 Finally, the reader may have noticed that there is a symmetry in the 
definition of a generalized path geometry, namely that the pair~
$\bigl(\cP,\cQ\bigr)$ defines a generalized path geometry on~$\Sigma$ if and
only if the pair~$\bigl(\cQ,\cP\bigr)$ also defines a generalized path 
geometry on~$\Sigma$.  This is known as the {\it dual\/} generalized path 
geometry.  It is not generally true that a generalized path geometry is
(locally) equivalent to its dual geometry.

\subhead 2.2.  The flat example: the projective plane \endsubhead
As a prelude to introducing Cartan's projective connection (a device for
solving the equivalence problem for generalized path geometries) in the 
next section, I will now discuss the path geometry of the flat example.

\subsubhead 2.2.1.  Lines in the projective plane \endsubsubhead
Let $V$ be $\bbR^3$, thought of as the space of column vectors of height~3
with real entries, and let $V^*$ be the space of row vectors of length~3 with
real entries.  The pairing $V^*\times V\to\bbR$ defined by matrix
multiplication is non-degenerate and will be taken to be the canonical 
pairing (thus justifying the notation~$V^*$). 

Let~$\bbS=\bigl(V\setminus\{0\}\bigr)/\bbR^+$ denote the space of 
{\it oriented\/} lines through the origin in~$V$.  This space is diffeomorphic
to the $2$@-sphere and canonically double covers~$\bbP^2 = 
\bigl(V\setminus\{0\}\bigr)/\bbR^*$, the space of (unoriented) lines through
the origin in~$V$.  If $v\in V$ is any non-zero vector, let~$[v]$
denote the corresponding point in~$\bbS$.  Note that, contrary to the usual
usage, $[-v]\not=[v]$.

Let~$\bbS^*=\bigl(V^*\setminus\{0\}\bigr)/\bbR^+$
denote the `dual' space of oriented lines through the origin in~$V^*$, with
$[\xi]$ denoting the oriented line corresponding to the non-zero element
~$\xi\in V^*$. 

Let~$\Sigma\subset\bbS^*\times\bbS$ denote the incidence correspondence
$$
\Sigma = \left\{ \>\bigl([\xi],[x]\bigr) \>\vert\>  \xi\cdot x = 0 \>\right\}.
$$
Then $\Sigma$ is a smooth $3$-manifold for which the two natural projections
$\lambda:\Sigma\to\bbS^*$ and $\pi:\Sigma\to\bbS$ are smooth submersions.
Moreover, $\Sigma$ carries a canonical generalized path geometry~
$\bigl(\cP,\cQ\bigr)$ where the leaves of the foliation~$\cP$ are the fibers
of~$\lambda$ and the leaves of the foliation~$\cQ$ are the fibers of~$\pi$.
The natural map~$\nu:\Sigma\to\bbP\bigl(T\bbS\bigr)$ is a $2$@-to@-$1$ 
covering, with the points~$\bigl([\pm\xi],[x]\bigr)$ going to the same 
point under the map~$\nu$.  Since the involution~$\bigl([\xi],[x]\bigr)\mapsto
\bigl([-\xi],[x]\bigr)$ maps the fibers of~$\lambda$ and~$\pi$
into fibers of these same maps, it follows that the generalized path
geometry defined on~$\Sigma$ descends to a well-defined (classical) path
geometry on~$\bbS$.  This is, of course, the classical geometry, where the
paths on~$\bbS$ are the `great circles'.

The group~$\SL(3,\bbR)$ acts on~$V$ on the left by matrix multiplication
and on~$V^*$ on the right by matrix multiplication.  This defines a left
action on~$V^*\times V$ by the rule
$$
g\cdot\bigl(\xi,x\bigr) = \bigl(\xi g^{-1},\>gx\bigr)
$$
for~$g\in\SL(3,\bbR)$.  This action preserves the locus $\xi\cdot x=0$ 
and descends to an action of~$\SL(3,\bbR)$ on~$\Sigma$ which commutes 
via the projections~$\lambda$ and $\pi$ with the natural actions of~
$\SL(3,\bbR)$ on~$\bbS^*$ and~$\bbS$ respectively.  Each of these actions is 
transitive and effective (i.e., only the identity in~$\SL(3,\bbR)$ acts as the 
identity transformation).  In particular, the group~$\SL(3,\bbR)$ acts
as a group of symmetries of the induced path geometry on~$\bbS$.  It is 
a classical fact that~$\SL(3,\bbR)$ is actually the full group of 
orientation preserving symmetries of this path geometry.

Let~$x_0={}^t\!(1,0,0)$ and let~$\xi_0=(0,0,1)$.  Let~$P_1\subset\SL(3,\bbR)$
be the stabilizer subgroup of~$[x_0]\in\bbS$, let $P_2\subset\SL(3,\bbR)$
be the stabilizer subgroup of~$[\xi_0]\in\bbS^*$, and let~$P=P_1\,\smcap\,P_2$.
Note that $P$ is the subgroup of upper triangular matrices in~$\SL(3,\bbR)$ 
with all diagonal entries positive.

Now define a map~$\tau:\SL(3,\bbR)\to\Sigma$ by the rule
$$
\tau(g) = g\cdot\bigl(\xi_0,x_0\bigr).
$$ 
The map~$\tau$ is surjective and its fibers are the left cosets of~$P$, 
so that~$\tau$ identifies~$\SL(3,\bbR)/P$ with~$\Sigma$.

Now let $\theta = g^{-1}\,dg$ be the canonical left-invariant 1@-form on~
$\SL(3,\bbR)$, i.e., $\theta$ is the unique $\eusl(3,\bbR)$@-valued 1@-form
on~$\SL(3,\bbR)$ whose value at the identity~$I_3\in\SL(3,\bbR)$ is 
the identity mapping
$$
\theta_{I_3}:T_{I_3}\SL(3,\bbR)=\eusl(3,\bbR)\to\eusl(3,\bbR).
$$
Since~$\eusl(3,\bbR)$ is the vector space of 3-by-3 matrices with trace
zero, $\theta$ can be expanded in the form
$$
\spreadmatrixlines{3pt}
\theta = 
\pmatrix 
\theta^0_0&\theta^0_1&\theta^0_2\\
\theta^1_0&\theta^1_1&\theta^1_2\\
\theta^2_0&\theta^2_1&\theta^2_2\\
\endpmatrix
$$
where~$\theta^0_0+\theta^1_1+\theta^2_2=0$, but the $\theta^a_b$ are
otherwise linearly independent.

 For comparison with the construction in the next section, I want to remark on
the following properties of~$\theta$ with respect to the oriented 
generalized path geometry on~$\Sigma$: 

 First, there is the fact that
$\SL(3,\bbR)$ acts on the left on all of the spaces~$\SL(3,\bbR)$, $\Sigma$,
$\bbS$, and $\bbS^*$ in such a way that it commutes with the maps $\tau$,
$\pi$, and $\lambda$; preserves the oriented generalized path geometry
on~$\Sigma$ and the orientations and path geometries on~$\bbS$ and $\bbS^*$;
and leaves~$\theta$ invariant. 

Second, if $\sigma:\Sigma\to\SL(3,\bbR)$ is any section 
of~$\tau:\SL(3,\bbR)\to\Sigma$ and $\phi = \sigma^*\theta$, 
then the components~$\phi^1_0$, $\phi^2_0$, and~
$\phi^2_1$ have the following properties:  First, the leaves of
the foliation~$\cP$ are the integral curves of~
$\phi^2_0=\phi^2_1=0$ while the 1@-form~
$\phi^1_0$ is $\cP$@-positive. 
Second, the leaves of the foliation~$\cQ$ are the integral curves of~
$\phi^1_0=\phi^2_0=0$ while the 1@-form~
$\phi^2_1$ is $\cQ$@-positive.

 Finally, note the Maurer-Cartan identity $d\theta = -\theta\w\theta$.

\subhead 2.3.  Cartan's projective connection \endsubhead
In \cite{Ca1}, Cartan introduced a device for determining when two generalized 
path geometries were equivalent, the so-called {\it projective connection\/}. 
This subsection contains an account of his results, modified slightly to take
into account the orientations. 

\subsubhead 2.3.1. The bundle and its connection form \endsubsubhead
Recall that the subgroup~$P\subset\SL(3,\bbR)$ has been defined to be 
the subgroup consisting of the upper triangular matrices with all diagonal 
entries positive.  Its Lie algebra is $\eup\subset\eusl(3,\bbR)$, the 
subspace consisting of upper triangular matrices with trace zero.

\proclaim{Theorem 2 (Cartan)}
Let $\bigl(\Sigma,\cP,\cQ\bigr)$ be an oriented generalized path geometry. 
There exists a principal right $P$@-bundle~$\tau:B\to\Sigma$ and an 
$\eusl(3,\bbR)$@-valued $1$@-form~$\theta$ on~$B$ 
with the following properties:
\roster
\item For each~$b\in B$, the map~$\theta_b:T_bB\to\eusl(3,\bbR)$ is
      an isomorphism and $\theta$ pulls back to each fiber of~$\tau$ to 
      be the canonical $\eup$@-valued left-invariant $1$@-form;
\item $R_g^*\theta = g^{-1}\theta g$ for each $g\in P$;
\item For some \rom(and hence any\rom) section~$\sigma:\Sigma\to B$, the 
      pullback $1$-form~$\phi = \sigma^*\theta$ has the properties 
      that, first, the leaves of the foliation~$\cP$ are the integral curves 
      of~$\phi^2_0=\phi^2_1=0$ while the $1$@-form~
      $\phi^1_0$ is $\cP$@-positive and, second, the leaves of the 
      foliation~$\cQ$ are the integral curves of~
      $\phi^1_0=\phi^2_0=0$ while the $1$@-form~
      $\phi^2_1$ is $\cQ$@-positive;
\item The curvature $2$@-form~$\Theta=d\theta+\theta\w\theta$ satisfies
        $$
        \Theta = 
        \pmatrix
        0&M\,\theta^1_0\w\theta^2_0&\Theta^0_2\\
        0&0& L\,\theta^2_1\w\theta^2_0\\
        0&0&0\\
        \endpmatrix
        $$
      for some functions~$L$ and $M$ on~$B$.
\endroster
The pair~$\bigl(B,\theta\bigr)$ is uniquely characterized by these
four properties: If $\bigl(B',\theta'\bigr)$ also satisfies them
then there exists a unique bundle isomorphism~$f:B\to B'$ covering the 
identity on~$\Sigma$ so that $f^*\theta'=\theta$. \qed
\endproclaim 

 For the proof, the reader can consult~\cite{Ca1} or the more modern 
treatment in~\cite{KN}. 

It should be noted that Cartan's proof is entirely constructive, i.e., 
he shows how to define~$B$ and~$\theta$ by an algorithmic process that only 
involves differentiation and solving certain explicit linear equations. 

In particular, suppose given a $3$@-manifold~$\Sigma$ on which 
there are defined two linearly independent 
vector fields~$X_1$ and $X_3$ with the property that $X_2 = [X_3,X_1]$ is 
linearly independent from $X_1$ and $X_3$. Consider the oriented generalized
path geometry~$\bigl(\Sigma,\cP,\cQ\bigr)$ defined by the oriented integral
curves of~$X_1$ and~$X_3$. Then the bundle~$B$ and the 1@-form~$\theta$
for the oriented generalized path geometry~$(\Sigma,\cP,\cQ)$ can be 
constructed by a (somewhat involved) recipe from the knowledge of the iterated
brackets of~$X_1$ and~$X_3$.  One does not need to be able to `integrate'
the vector fields~$X_1$ or~$X_3$.

\subsubhead 2.3.2. Some applications \endsubsubhead
 From Theorem~2, Cartan draws several conclusions. 

 First, two oriented
generalized path geometries $\bigl(\Sigma_1,\cP_1,\cQ_1\bigr)$  and
$\bigl(\Sigma_2,\cP_2,\cQ_2\bigr)$ are equivalent if and only if there exists
a diffeomorphism~$f:B_1\to B_2$ satisfying $f^*\theta_2 = \theta_1$. 

Second, the group of symmetries of a given oriented generalized path geometry~
$\bigl(\Sigma,\cP,\cQ\bigr)$ is isomorphic to the group of diffeomorphisms
$f:B\to B$ that satisfy~$f^*\theta=\theta$.  Since the components of~$\theta$
provide a coframing of~$B$, it follows (see~\cite{Ko}) that this 
diffeomorphism group has the structure of a Lie group of dimension at most 8.

Third, for the flat example, $\Sigma\subset\bbS^*\times\bbS$, the pair~
$\bigl(B,\theta\bigr)$ can be taken to be $\bigl(\SL(3,\bbR),g^{-1}\,dg\bigr)$
while~$\Sigma$ is identified with~$\SL(3,\bbR)/P$.  In particular, in
this case, since~$\Theta=d\theta+\theta\w\theta=0$, it follows that the
functions~$L$~and~$M$ vanish identically on the flat example.

Thus, if $\bigl(\Sigma,\cP,\cQ\bigr)$ is to be locally equivalent to the flat 
example, the functions $L$ and $M$ must vanish identically.
Conversely, if $L=M=0$, then the Bianchi identity $d\Theta = \Theta\w\theta
-\theta\w\Theta$ implies that $\Theta^0_2$ vanishes identically as well, so
that~$\Theta$ itself vanishes identically, i.e., $d\theta = -\theta\w\theta$.
It follows that if~$\Sigma$ is simply connected, then there exists a 
smooth mapping~$f:B\to\SL(3,\bbR)$ which satisfies $f(b\cdot g)=f(b)\cdot g$ 
for all $g\in P$ and which pulls back the canonical left-invariant form on~
$\SL(3,\bbR)$ to~$\theta$ on~$B$. This mapping is a local diffeomorphism and 
is unique up to left translation by a constant in~$\SL(3,\bbR)$. In particular,
this induces a local diffeomorphism~$\bar f:\Sigma\to \SL(3,\bbR)/P$ which is 
a local equivalence of oriented generalized path geometries.  Thus, the
vanishing of~$L$ and~$M$ is also a sufficient condition that 
$\bigl(\Sigma,\cP,\cQ\bigr)$ be locally equivalent to the flat example.

 Fourth, the Bianchi identity shows that $L$ and $M$ cannot be constant unless 
they are identically zero.  Since any diffeomorphism~$f:B\to B$ that fixes~
$\theta$ must also fix~$\Theta$ and hence~$L$ and $M$, it follows that, 
if $L$ or $M$ is non-zero, then the group of symmetries of the generalized 
path geometry has dimension strictly less than 8.

\subsubhead 2.3.3. Computations \endsubsubhead
Let~$\bigl(\Sigma,\cP,\cQ\bigr)$ be an oriented generalized path geometry
with Cartan structure bundle and connection~$(B,\theta)$.
If $\sigma:\Sigma\to B$ is any section, then~$\phi=\sigma^*\theta$
is an $\eusl(3,\bbR)$@-valued $1$@-form with four properties:
\roster
\item  $\phi^1_0\w\phi^2_0\w\phi^2_0$
       is a non-vanishing 3@-form on~$\Sigma$; 
\item  the leaves of the foliation~$\cP$ are the integral curves 
       of~$\phi^2_0=\phi^2_1=0$ while the 1@-form~$\phi^1_0$ 
       pulls back to each such leaf to be a $\cP$@-positive 1@-form; 
\item  the leaves of the foliation~$\cQ$ are the integral curves 
       of~$\phi^1_0=\phi^2_0=0$ while the 1@-form~$\phi^2_1$ 
       pulls back to each such leaf to be a $\cQ$@-positive 1@-form;
\item  the curvature $2$@-form~$\Phi=d\phi+\phi\w\phi$ satisfies
       $$
       \spreadmatrixlines{3pt}
       \Phi = 
       \pmatrix
       0&\overline{M}\,\phi^1_0\w\phi^2_0&\Phi^0_2\\
       0&0& \overline{L}\,\phi^2_1\w\phi^2_0\\
       0&0&0\\
       \endpmatrix
       $$
       for some functions~$\overline{L}$ and $\overline{M}$ on~$\Sigma$.
\endroster

Conversely, if~$\phi$ is any $\eusl(3,\bbR)$@-valued $1$@-form on~$\Sigma$ 
with these four properties, then it is of the form~$\sigma^*\theta$
for some section~$\sigma$ of~$B$.

In particular, in order to show that a given oriented generalized path
geometry is locally equivalent to the flat example, it suffices to 
construct a 1@-form~$\phi$ on~$\Sigma$ satisfying the first three
properties and
$$
\Phi 
= d\phi + \phi\w\phi
=
\pmatrix
0&0&\Phi^0_2\\
0&0&0\\
0&0&0\\
\endpmatrix\,.
$$
The Bianchi identity $d\Phi = \Phi\w\phi-\phi\w\Phi$ then implies
that the component~$\Phi^0_2$ must also be zero.

\subhead 2.4.  The generalized path geometry of a Finsler structure \endsubhead
To each generalized Finsler structure~$\bigl(\Sigma,\omega\bigr)$ there
is canonically associated an oriented generalized path geometry~
$\bigl(\Sigma,\cP,\cQ\bigr)$.  The leaves of~$\cP$ are 
the integral curves of~$\Xbold_1$ and the leaves of~$\cQ$ are the integral 
curves of~$\Xbold_3$.  An orientation is fixed by declaring that~$\omega_1$ 
be $\cP$@-positive and $\omega_3$ be $\cQ$@-positive.

\subsubhead 2.4.1. The projective connection form \endsubsubhead
Keeping the structure equations and Bianchi identity notation as given 
in~\S1.3, consider the $\eusl(3,\bbR)$@-valued $1$@-form
$$
\spreadmatrixlines{5pt}
\phi = 
\pmatrix
{\ts{1\over3}}\bigl(I\,\omega_3-J\,\omega_2\bigr)
&-K\,\omega_1-{\ts{1\over3}}K_3\,\omega_2
&-{\ts{2\over3}}K_3\,\omega_1-U\,\omega_2+{\ts{1\over3}}(I_2{+}J_3)\,\omega_3\\
\omega_1&{\ts{1\over3}}\bigl(I\,\omega_3-J\,\omega_2\bigr)
&-{\ts{1\over3}}(I_2{+}J_3)\,\omega_2-\omega_3\\
\omega_2&\omega_3&-{\ts{2\over3}}\bigl(I\,\omega_3-J\,\omega_2\bigr)\\
\endpmatrix
$$
where~$U=K+{1\over3}K_{33}+{1\over3}K_3\,I$.  A short calculation yields
$$
d\phi+\phi\w\phi
=
\pmatrix
0&\overline{M}\,\omega_1\w\omega_2&\Phi^0_2\\
0&0&\overline{L}\,\omega_3\w\omega_2\\
0&0&0\\
\endpmatrix\,,
$$
where
$$
\align
3\overline{M} &= -K_{31}+3K_2\,,\\
3\overline{L} &= -I_{23}-J_{33}-2I(I_2+J_3)-6J\,.\\
\endalign
$$
Since $\phi$ has all four properties listed in~\S2.3.3, it is
the pullback of the projective connection form via a section of the
projective structure bundle.

\subsubhead 2.4.2. Projective flatness\endsubsubhead
The results of the last subsection coupled with the discussion in~
\S\S2.3.2-3 now combine to show that the generalized path geometry 
associated to a generalized Finsler structure is locally equivalent
to that of the flat example if and only if its invariants satisfy
$$
\align
                 K_{31}-3K_2 &=0\,,\\
I_{23}+J_{33}+2I(I_2+J_3)+6J &=0\,.\\
\endalign
$$
This characterization goes back at least to Rund~\cite{Ru, p.~261}, but in 
some form is due to Berwald~\cite{Be}, at least in the context of 
classical Finsler structures.  (Since the characterization is local, the 
distinction between the classical and the generalized case is not important 
here.)

In the context of Hilbert's Fourth Problem, this gives a diffeomorphism
invariant characterization of the Finsler Lagrangians whose geodesics
can be mapped to the straight lines by some diffeomorphism. 
A generalized Finsler structure that satisfies these conditions will
be said to be {\it projectively flat\/}.

\head3. Projectively flat Finsler structures on~$S^2$ satisfying $K=1$\endhead

The main object of this article is to classify the Finsler structures 
on~$S^2$ that satisfy~$K=1$ and are projectively flat.  However, it
is useful to start in a more general setting, as it leads to a stronger
result. 

Thus, throughout this section, $\bigl(\Sigma,\omega\bigr)$ will 
be a connected, compact $3$@-manifold~$\Sigma$ with finite fundamental
group endowed with a generalized Finsler structure~$\omega$ satisfying the 
conditions that $K\equiv1$ and that the induced generalized path geometry 
on~$\Sigma$ be projectively flat. 
The case where~$\Sigma$ is an actual Finsler structure on~$S^2$ that 
satisfies~$K\equiv1$ and is projectively flat is an example, for such a
$\Sigma$ will be compact and its fundamental group will be~$\bbZ_2$. 
 
It is not clear that the assumption that~$\Sigma$ have finite 
fundamental group is necessary for the conclusions of the theorems, 
as it is possible that compactness implies it.  However, this assumption 
does simplify the discussion and is satisfied in the cases of interest, 
so it seems reasonable to include it.

It may seem that projective flatness plus constant curvature is so 
stringent an assumption that there will not be many such structures, even
locally.  However, surprisingly, it turns out that there are very many
local examples, as was discovered by Funk~\cite{Fu1,Fu2}.  In fact, the 
problem of characterizing the global solutions remained a concern throughout 
 Funk's career~\cite{Fu3}. 

In this section, the first result will be Theorem~3, which asserts
that any generalized Finsler structure~$\bigl(\Sigma,\omega\bigr)$ on
a simply connected, compact $3$@-manifold that is projectively flat
and satisfies~$K=1$ is the double cover of a classical Finsler structure
on~$\bbS$ with these properties.  From this, the classification in
the case of finite fundamental group follows.

The second result is a finiteness result, Theorem~4, which shows that,
in a certain sense, two such structures that agree to sufficiently high
order at corresponding points must be globally equivalent.

The final result is an existence result, Theorem~5, showing that
non-Riemannian examples exist and depend on two essential parameters.

One surprising aspect of the proof is that it shows that the compact examples
are naturally an open subset of an analytic family of examples, not all
of which are compact.

\subhead 3.1. The structure equations \endsubhead
Using the condition~$K=1$, the structure equations for the canonical
coframing~$\omega$ on~$\Sigma$ simplify to
$$
\aligned
d\omega_1 &= -\omega_2\w\omega_3\,,\\
d\omega_2 &= \phantom{-}\bigl(\omega_1 - I\,\omega_2\bigr)\w\omega_3\,,\\
d\omega_3 &=           -\bigl(\omega_1 - J\,\omega_3\bigr)\w\omega_2\,,\\
\endaligned 
\tag1
$$
and the Bianchi identities reduce to:
$$
\alignedat4
dI ={}&&  J\,\omega_1 &&{}+ I_2\,\omega_2 &&{}+ I_3\,\omega_3&{}\,,\\
dJ ={}&& -I\,\omega_1 &&{}+ J_2\,\omega_2 &&{}+ J_3\,\omega_3&{}\,.\\
\endaligned
\tag2
$$
Since~$K\equiv1$, the first of the two conditions for projective flatness
of the associated generalized path geometry, namely, $K_{31}-3K_2 =0$,
is an identity, so the only remaining condition to impose is
$I_{23}+J_{33}+2I(I_2+J_3)+6J =0$, so I assume this from now on. 

The quantity~${1\over3}(I_2{+}J_3)$ turns up many places in what follows, 
so I will denote it by~$T$.  The condition for projective flatness then 
becomes~
$$
T_3 + 2IT + 2 J = 0.
\tag3
$$ 

Now, the structure equations yield~$d(I\,\omega_3-J\,\omega_2)
=3T\,\omega_2\w\omega_3$.  Applying the exterior derivative to both
sides of this equation and using the structure equations again yields
$0 = 3\,dT\w\omega_2\w\omega_3$, so it follows that $T_1=0$.  Thus,
$$
dT = T_2\,\omega_2 - 2(IT+J)\,\omega_3\,.
$$
Differentiating this relation and wedging both sides with $\omega_2$ yields
$$
0 = \bigl(T_2\,\omega_1\w\omega_3-2(JT-I)\,\omega_1\w\omega_3\bigr)
\w\omega_2\,,
$$
so that $T_2 = 2(JT-I)$, i.e.,
$$
dT = 2(JT-I)\,\omega_2 - 2(IT+J)\,\omega_3\,.
\tag4
$$

\subhead 3.2. The canonical immersion \endsubhead
Now consider the $\eusl(3,\bbR)$@-valued $1$@-form
$$
\spreadmatrixlines{5pt}
\phi = 
\pmatrix
{\ts{1\over3}}\bigl(I\,\omega_3-J\,\omega_2\bigr)
&-\omega_1
&-\omega_2+T\,\omega_3\\
\omega_1&{\ts{1\over3}}\bigl(I\,\omega_3-J\,\omega_2\bigr)
&-\omega_3-T\,\omega_2\\
\omega_2&\omega_3&-{\ts{2\over3}}\bigl(I\,\omega_3-J\,\omega_2\bigr)\\
\endpmatrix\,.
\tag5
$$
By equations~(1--4), the $\eusl(3,\bbR)$@-valued $1$@-form~$\phi$ satisfies~
$d\phi=-\phi\w\phi$. 

\proclaim{Theorem 3}
Suppose that~$\Sigma$ is a compact, simply connected $3$@-manifold
on which there is a generalized Finsler structure~$\omega$ that satisfies
$K=1$ and is projectively flat.  Define~$\phi$ as above. 
\roster
\item 
There exists a smooth immersion~$g:\Sigma\to\SL(3,\bbR)$, unique
up to left translation by a constant matrix in~$\SL(3,\bbR)$, that satisfies
$\phi = g^{-1}\,dg$. 
\item
The map~$g$ is a double cover onto its image~$\bar\Sigma=g\bigl(\Sigma\bigr)$, 
a smoothly embedded $3$@-dimensional submanifold of~$\SL(3,\bbR)$. 
\item
The quotient generalized Finsler structure on~$\bar\Sigma\simeq\Sigma/\bbZ_2$ 
is both amenable and geodesically amenable, with the projection~
$\bar\Sigma\to\bbS=\SL(3,\bbR)/P_1$, resp.~$\bar\Sigma\to\bbS^* = 
\SL(3,\bbR)/P_2$, having the leaves of the induced foliation~$\bar\cQ$, 
resp.~$\bar\cP$, as fibers. 
\item
The induced canonical immersion~$\bar\nu:\bar\Sigma\to T\bbS$ embeds 
$\bar\Sigma$ as a Finsler structure on~$\bbS$ whose geodesics are the 
standard lines in~$\bbS$.
\endroster
\endproclaim

\demo{Proof}  Since $\Sigma$ is simply connected and $\phi$ satisfies
the Maurer-Cartan equation $d\phi+\phi\w\phi=0$ (sometimes known as the
`zero-curvature equation'), it follows by Cartan's technique of the graph
that there exists a smooth 
mapping~$g:\Sigma\to\SL(3,\bbR)$ so that ~$\phi = g^{-1}\,dg$ 
and that this mapping is unique up to left translation, i.e., if $\phi 
= h^{-1}\,dh$ for some other map~$h:\Sigma\to\SL(3,\bbR)$ then there 
exists a constant matrix~$h_0\in\SL(3,\bbR)$ so that~$h = h_0\,g$.

Since the three subdiagonal entries of~$g^{-1}\,dg=\phi$
form a coframing of~$\Sigma$, it follows that~$g$ is an immersion. 

Now, recall that $\bbS = \SL(3,\bbR)/P_1$ where the mapping~
$[x] = \pi\circ\tau:\SL(3,\bbR)\to\bbS$ defined by~$[x](g) = [g\cdot x_0]$
has the left cosets of~$P_1$ as fibers.  From the definition of~$\phi$, it
then follows that~$[x]\circ g:\Sigma\to\bbS$ is a submersion whose
fibers are unions of integral curves of~$\Xbold_3$, i.e., the leaves of~$\cQ$.

Since~$\Sigma$ is compact, $[x]\circ g$ must be a 
surjective submersion and each fiber of~$[x]\circ g$ must be compact.
In particular, each such fiber must be diffeomorphic to a disjoint union
of a finite number of circles.  Since each leaf of~$\cQ$ must be a component
of such a fiber, each leaf of~$\cQ$ must be compact.

The submersion $[x]\circ g:\Sigma\to\bbS$ must therefore be a fibration. 
The exact homotopy sequence of this fibration includes the segment
$$
\{{\bold1}\}\simeq\pi_1\bigl(\bbS\bigr)\longrightarrow
\pi_0\bigl(F\bigr)\longrightarrow
\pi_0\bigl(\Sigma\bigr)\simeq\{{\bold1}\}\,,
$$
where $F$ is a typical fiber.  Thus $\pi_0\bigl(F\bigr)$ must be trivial, 
i.e., the fibers of~$[x]\circ g$ are connected, so that each fiber consists 
of precisely one leaf of~$\cQ$.

In sum, the generalized Finsler structure on~$\Sigma$ is amenable
and I can regard the submersion $[x]\circ g:\Sigma\to\bbS$ as the 
canonical $\cQ$@-quotient mapping. (By the way, it is {\it not\/} 
generally true that the integral curves of~$\Xbold_3$ must be amenable
for a generalized Finsler structure satisfying only~$K=1$, even when the 
underlying manifold~$\Sigma$ is compact with finite fundamental group. 
The hypothesis of projective flatness has been used in a crucial way in 
this argument.)
 
Let~$\nu:\Sigma\to T\bbS$ be the induced canonical immersion as in~\S1.2. 
By Proposition~1, it follows that the image of
each fiber~$\Sigma_s = ([x]\circ g)^{-1}(s)$ in~$T_s\bbS$ is a closed curve 
that is locally strictly convex towards~$0_s$ for all~$s\in\bbS$. 
The number of times that~$\nu\bigl(\Sigma_s\bigr)$ winds around~$0_s$ is
independent of~$s$.  Moreover, this winding number must be~2 since the 
fundamental group of the unit circle bundle~$U\bbS\subset T\bbS$ of any
Riemannian metric~$g$ on~$\bbS$ is~$\bbZ_2$ and the natural mapping~
$\beta:\Sigma\to U\bbS$ defined by~$\beta(\ubold)=\nu(\ubold)/|\nu(\ubold)|_g$ 
is a covering covering map and hence must be $2$@-to@-$1$.

Now consider the dual map~$[\xi]\circ g:\Sigma\to\bbS^*=\SL(3,\bbR)/P_2$.
This mapping is also a submersion and its fibers are unions of integral curves 
of~$\Xbold_1$. The same argument as before applies to show that~$[\xi]\circ g$ 
must be a surjective submersion  with each of its fibers being connected and 
hence consisting of a single, closed integral curve of~$\Xbold_1$.  Thus, the 
generalized Finsler structure on~$\Sigma$ is also geodesically amenable. 

Now, the leaves of $\cP$ satisfy $\omega_2=\omega_3=0$ and hence, from the 
definition of~$\phi$, it follows that each leaf of~$\cP$ maps via~$g$ 
submersively onto a left coset of the circle subgroup
$$
H = \left\{
\pmatrix
\cos\theta&-\sin\theta&0\\
\sin\theta&\phantom{-}\cos\theta&0\\
0&0&1\\
\endpmatrix
\>\vrule\>
\theta\in\bbR\>
\right\}.
$$
The image of such a left coset under the quotient map~$[x]:\SL(3,\bbR)\to
\bbS$ is a projective `line', i.e., the set of oriented lines lying
in a plane in~$V$.  In particular, it follows that the geodesics of the
generalized Finsler structure on~$\Sigma$ map to the standard lines on~$\bbS$.

Let~$\Phi_t:\Sigma\to\Sigma$ be the time~$t$ flow of the vector 
field~$\Xbold_1$.  Since~$\Sigma$ is compact, this flow exists for all
time.  Since the structure equations show that $\Lie_{\Xbold_1}\omega_1=0$, 
$\Lie_{\Xbold_1}\omega_2=\omega_3$, and $\Lie_{\Xbold_1}\omega_3=-\omega_2$, 
it follows that
$$
\align
\Phi^*_s\omega_1 & = \omega_1\,,\\
\Phi^*_s\omega_2 & = \phantom{-}\cos s\,\omega_2 + \sin s\,\omega_3\,,\\
\Phi^*_s\omega_3 & =           -\sin s\,\omega_2 + \cos s\,\omega_3\,.\\
\endalign
$$ 
In particular, $\Phi_{2\pi}:\Sigma\to\Sigma$ satisfies~
$\Phi^*_{2\pi}\omega_i=\omega_i$, for~$i=1,2,3$.  Consequently,
$\Phi^*_{2\pi}$ must also preserve~$I$, $J$, $T$, and hence~$\phi$. It is
also evident that
$$
g\bigl(\Phi_s(\ubold)\bigr) = g(\ubold)\cdot
\pmatrix\cos s&-\sin s&0\\\sin s&\phantom{-}\cos s&0\\0&0&1\\\endpmatrix,
$$
so $g\bigl(\Phi_{2\pi}(\ubold)\bigr) = g(\ubold)$ for all~$\ubold\in\Sigma$.
In particular, $\nu\circ\Phi_{2\pi}=\nu$, so that~$\Phi_{2\pi}$ acts
as a deck transformation for the double cover~$\beta:\Sigma\to U\bbS$
constructed above.  It follows either that~$\Phi_{2\pi}$ is the identity
or else that it is a fixed point free involution.

Now~$\Phi_{2\pi}$ cannot be the identity, since the image of a left coset 
of~$H$ is a line in~$\bbS$ whose tangential lift to $U\bbS$ is a generator 
of its fundamental group~$\bbZ_2$.  Thus the double cover~$\Sigma\to U\bbS$ 
must be a non-trivial double cover over each such lift, implying that the 
map~$g$ is $2$@-to@-$1$ on each leaf of~$\cP$.

In sum, the map~$g$ commutes with the fixed point free involution~$\Phi_{2\pi}$ 
and hence induces a smooth immersion of the quotient~$\bar\Sigma=\Sigma/\bbZ_2$
into~$\SL(3,\bbR)$.  This immersion must be an embedding since the induced map
$\bar\nu:\bar\Sigma\to U\bbS$ is now seen to be a $1$@-sheeted covering and
factors through~$g$. \qed
\enddemo 

\proclaim{Corollary 1}
If $\bigl(\Sigma,\omega\bigr)$ is a generalized Finsler structure on a
connected compact manifold~$\Sigma$ with finite fundamental group that 
satisfies $K\equiv1$ and whose associated path geometry is locally flat, 
then either $\bigl(\Sigma,\omega\bigr)$ is a quotient of a projectively
flat Finsler structure on~$\bbS$ by a finite subgroup of Finsler isometries
or else it is a double cover of such a quotient.
\endproclaim

\demo{Proof}
Let~$\tilde\Sigma\to\Sigma$ be the simply connected cover and let~$\tilde
\omega$ be the pull-back of the generalized Finsler structure~$\omega$
via the covering map.  Since the deck transformations of this covering map
preserve~$\tilde\omega$, they must preserve the dual vector fields and,
in particular, they must commute with the involution~$\Phi_{2\pi}:
\tilde\Sigma\to\tilde\Sigma$ constructed in the proof of
the above theorem.  There are now two cases to consider:

 First, if $\Phi_{2\pi}$ is an element of the deck transformation group for
$\tilde\Sigma\to\Sigma$, then it lies in the center of this group and
so covers the identity on~$\Sigma$.  Dividing by the action of~$\Phi_{2\pi}$
then induces a covering map~$\bar\Sigma\to\Sigma$, so that~$\Sigma$ is
a quotient of the Finsler structure~$\bar\Sigma$ by the induced action
of a finite group of Finsler isometries acting on~$\bbS$.

Second, if $\Phi_{2\pi}$ is not an element of the deck transformation group,
then because it commutes with the action of this group on~$\tilde\Sigma$,
it follows that there is a non-trivial involution~$F:\Sigma\to\Sigma$
which preserves~$\omega$ and makes the following diagram commute:
$$
\CD
\tilde\Sigma   @>{\Phi_{2\pi}}>> \tilde\Sigma \\
@VVV @VVV\\
\Sigma @>F>> \Sigma\\
\endCD
$$
Since~$F$ preserves~$\omega$ and is not the identity, it cannot have
any fixed points.

It follows that there exists a quotient~$\Sigma\to \hat\Sigma$ whose
deck transformation group is generated by the involution~$F$.  A simple
diagram chase now shows that $\hat\Sigma$ is a quotient of~$\bar\Sigma
=\tilde\Sigma/\bbZ_2$. \qed
\enddemo
 
\subhead 3.3. The induced Riemannian metric on~$\bbS^*$ \endsubhead
In view of this theorem and its corollary, no essential generality is lost
by assuming that~$\bigl(\Sigma,\omega\bigr)$ is a Finsler structure 
on~$\bbS$ whose geodesics are the standard ones and that~$\Sigma$ has
been embedded into~$\SL(3,\bbR)$ via a mapping~$g:\Sigma\to\SL(3,\bbR)$
that satisfies~$(5)$ so I will be assuming this from now on.

Since the integral curves of~$\Xbold_1$ are periodic of period~$2\pi$
on~$\Sigma$, it follows that the flow of~$\Xbold_1$ defines a free action
of the unit circle on~$\Sigma$ whose orbits are the leaves of~$\cP$, i.e.,
the fibers of the submersion~$[\xi]:\Sigma\to\bbS^*$.  Thus, $\Sigma$ can
be regarded as a principal $S^1$@-bundle over~$\bbS^*$.

The structure equations will now show that there are certain semi-basic
tensors for this projection which are invariant under this $S^1$@-action
and hence yield well-defined tensors on~$\bbS^*$.
 
 For example, since $\Lie_{\Xbold_1}\omega_2=\omega_3$ and 
$\Lie_{\Xbold_1}\omega_3=-\omega_2$, it follows that 
$$
\Lie_{\Xbold_1}\,\omega_3\w\omega_2 = 0
\qquad\text{and}\qquad
\Lie_{\Xbold_1}\,\bigl((\omega_3)^2+(\omega_2)^2\bigr) = 0.
$$
In particular, there exists a $2$@-form~$dA$ on~$\bbS^*$ so that 
$[\xi]^*(dA) = \omega_3\w\omega_2$ and there exists a Riemannian
metric~$ds^2$ on~$\bbS^*$ so that $[\xi]^*(ds^2) = (\omega_3)^2+(\omega_2)^2$.
 From now on, I will consider~$\bbS^*$ to be an oriented Riemannian manifold
with metric~$ds^2$ and area form~$dA>0$. 

Now the structure equations show that the map~$\mu:\Sigma\to T\bbS^*$
defined by the formula~$\mu(\ubold)=[\xi]'(\ubold)\bigl(\Xbold_3(\ubold)\bigr)$
embeds~$\Sigma$ as the unit circle bundle of~$\bbS^*$ endowed with the
metric~$ds^2$.  In particular, the $1$@-forms~$\omega_3$ and $\omega_2$
can be regarded as the tautological $1$@-forms on this circle bundle.  The
structure equations then show that the $1$@-form
$$
\rho = -\omega_1 + I\,\omega_2 + J\,\omega_3
$$
is the Levi-Civita connection for~$\Sigma$ regarded as this unit circle bundle.
Note that this implies that the Gauss curvature~$K$ of this metric, which, by
definition, satisfies~$d\rho= -K\,\omega_3\w\omega_2$, must be
$$
K = 1 - I_3 + J_2-I^2-J^2\,.
$$

\subhead 3.4. A complex form of the structure equations \endsubhead
Since~$\bbS^*$ inherits a canonical metric and orientation from~$\Sigma$,
it follows that it has a unique structure as a Riemann surface for which
the metric~$ds^2$ is conformal and the area form~$dA$ is a positive
$(1,1)$@-form.  A complex valued 1@-form~$\alpha$ on~$\bbS^*$ is of
type~$(1,0)$ with respect to this complex structure if and only if it
satisfies~$[\xi]^*(\alpha) = A\,(\omega_3+i\,\omega_2)$ for some function~
$A$ on~$\Sigma$. 

 For simplicity, set~$\zeta = \omega_3+i\,\omega_2$. 
Then the structure equations imply
$$
d\zeta = -i\rho\w\zeta
\tag6
$$
where, as before, $\rho = -\omega_1{+}I\,\omega_2{+}J\,\omega_3$. 

In this notation, equation~\S3.1.4 becomes
$$
dT = -(T-i)(I+iJ)\zeta-(T+i)(I-iJ)\overline\zeta.
\tag7
$$
Taking the exterior derivative of both sides of this equation and 
using the structure equations then yields
$$
0 = 2\bigl(J_2-I_3 + 3 (T^2{-}I^2{-}J^2)\bigr)\,\omega_3\w\omega_2\,,
$$
so it follows that $I_3-J_2 = 3(T^2-I^2-J^2)$.   Using this relation, 
the equations~$(2)$ for $dI$ and $dJ$ can be written in the form
$$
d\left(\,(I+iJ)\,\zeta\,\right) = 
{\ts{1\over2}}\bigl(2I^2+2J^2-3T^2-3iT\bigr)\,\zeta\w\overline\zeta.
\tag8
$$

Note also that $K = 1 + 2I^2 + 2J^2-3T^2$, so that,
$$
d\rho = -{\ts{i\over2}}\bigl(1+2I^2+2J^2-3T^2\bigr)\,\zeta\w\overline\zeta.
\tag9
$$

Conversely, the satisfaction of equations (6--9) is equivalent to the
satisfaction of the system~$d\phi+\phi\w\phi=0$, once the translation of
notation is made.

\subhead 3.5. A holomorphic cubic differential \endsubhead
I now want to investigate the consequences of equations~(6--9). 
To this end, it turns out to be useful to define $a = -(I+iJ)/(T+i)$.
(Since~$T$ is real, $T{+}i$ can never vanish, so this formula
does indeed define~$a$ smoothly.)  In terms of~$a$, the structure equations
now take the simpler form
$$
\aligned
d\zeta &= -i\rho\w\zeta,\\
d\rho&=-{\ts{i\over2}}(1+2(T^2+1)|a|^2-3T^2)\,\zeta\w\overline\zeta\,\\
dT &= (T^2+1)\bigl(a\,\zeta+\overline a\,\overline\zeta\bigr)\,\\ 
d\bigl(\,a\,\zeta\,\bigr) &= {\ts{3\over2}}T\,\zeta\w\overline\zeta\,.\\
\endaligned
\tag10
$$
Expanding and rearranging this last equation yields
$$
\bigl(\,da - ia\,\rho + {\ts{3\over2}}T\,\overline\zeta\,\bigr)\w\zeta = 0.
$$
It follows that there exists a complex valued function~$b$ on~$\Sigma$ so
that
$$
da = ia\,\rho +(b+a^2T)\,\zeta - {\ts{3\over2}}T\,\overline\zeta.
\tag11
$$
(Naming the $\zeta$@-coefficient $(b{+}a^2T)$ instead of~$b$ 
simplifies the next formula.)  Differentiating (11) and using the 
structure equations yields
$$
0 = \bigl(\,db - 2ib\,\rho + a\,\overline\zeta\,\bigr)\w\zeta.
$$
Thus, there exists a complex valued function~$c$ on~$\Sigma$ so that
$$
db = 2ib\,\rho +(c+2abT-{\ts{2\over9}}a^3)\,\zeta - a\,\overline\zeta.
\tag12
$$
(Naming the $\zeta$@-coefficient in this equation~
$(c{+}2abT{-}{\ts{2\over9}}a^3)$ instead of~$c$ simplifies the next 
formula).  Finally, differentiating (12) and using
the equations found so far yields the decisive formula
$$
0 = \bigl(\,dc - 3ic\,\rho\,\bigr)\w\zeta.
\tag13
$$
 From this it follows that there is a cubic differential~$\gamma$ 
on~$\bbS^*$ satisfying $[\xi]^*(\gamma) = c\,\zeta^3$ and, moreover, that
this differential form is holomorphic.

\subhead 3.6. A vanishing theorem and it consequences \endsubhead
Now~$\bbS^*$ is topologically a $2$@-sphere, so by the Uniformization Theorem 
it must be equivalent as a Riemann surface to the Riemann sphere~$\bbC\bbP^1$.
Since there are no non-zero holomorphic differentials
of positive degree on~$\bbC\bbP^1$, it follows that~$\gamma\equiv0$, or
equivalently, $c\equiv0$.  Thus, in complex form, the structure equations 
on~$\Sigma$ have finally been reduced to
$$
\aligned
d\zeta &= -i\rho\w\zeta\,,\\
d\rho&=-{\ts{i\over2}}(1+2(T^2+1)|a|^2-3T^2)\,\zeta\w\overline\zeta\,,\\
dT &= (T^2+1)\bigl(a\,\zeta+\overline a\,\overline\zeta\bigr)\,,\\ 
da &=  ia\,\rho +(b+a^2T)\,\zeta - {\ts{3\over2}}T\,\overline\zeta\,,\\
db &= 2ib\,\rho +(2abT-{\ts{2\over9}}a^3)\,\zeta - a\,\overline\zeta\,.\\
\endaligned
\tag14
$$
Note that these formulae give expressions for the exterior 
derivative of every quantity appearing in the formulae.  Moreover, exterior
differentiation of these formulae simply yield identities, so that there
are no further relations to be found by this method.

I will refer to the map~$(T,a,b):\Sigma\to \bbR\times\bbC\times\bbC$
as the {\it structure function\/} of the generalized Finsler structure~
$\omega$. 

In the next three subsections, I am going to prove a uniqueness theorem
and an existence theorem.  Since the constructions are somewhat long,
I am going to state the results here, hoping that knowing the results
will provide motivation to the reader for following the constructions.

\proclaim{Theorem 4}
Suppose that~$\Sigma$ and $\Sigma^*$ are two connected, simply-connected
and compact $3$-manifolds endowed with generalized Finsler structures~
$\omega$ and $\omega^*$, respectively, that are each projectively flat
and satisfy~$K=1$.  Let~$\bigl(T,a,b\bigr)$ and $\bigl(T^*,a^*,b^*\bigr)$
respectively denote the real and complex-valued functions constructed 
from~$\omega$ and $\omega^*$ by the formulae in \S$3.4$--$5$.
If there exist points~$\ubold\in\Sigma$ and $\ubold^*\in\Sigma^*$ for
which the equalities
$$
T(\ubold)=T^*(\ubold^*),\qquad a(\ubold)=a^*(\ubold^*),\qquad
b(\ubold)=b^*(\ubold^*)
$$ 
all hold, then there exists a unique diffeomorphism $f:\Sigma^*\to\Sigma$
so that~$f^*(\omega) = \omega^*$ and $f(\ubold^*) = \ubold$.
\endproclaim 

This uniqueness is accompanied by the following existence result:

\proclaim{Theorem 5}
There exists an open neighborhood~$U$ of~$(0,0,0)$ 
in~$\bbR\times\bbC\times\bbC$ so that for any~$\bigl(T_0,a_0,b_0\bigr)
\in U$ there exists a generalized Finsler structure~$\omega$ on~$S^3$
that is projectively flat, satisfies~$K=1$, and whose structure function
~$(T,a,b)$ assumes the value~$\bigl(T_0,a_0,b_0\bigr)$.
\endproclaim

Now, from the equations~(14), the rank of the differential of the structure 
function can be computed in terms of the value of the structure function
itself.  Indeed, calculation shows that the structure function is an
immersion except at the points where either $a=0$ or else where there is
a real number~$r$ so that $b = a^2r$ and $|a|^2 = {9\over2}(1-3rT)/(1+9r^2)$.
Since the generic point~$(T_0,a_0,b_0)$ in $\bbR\times\bbC\times\bbC$ does
not satisfy either of these conditions, Theorem~5 implies
that for the `generic' generalized Finsler structures satisfying equations~
(14), the rank of the structure functions is 3 at most places.  In other words,
the images of the various structure functions are of dimension 3 at most
places.  This suggests that the images of these structure functions should
be the simultaneous level sets of a pair of functions defined on~
$\bbR\times\bbC\times\bbC$. 

This is, indeed, the case.  A little experimentation with equations~(14) 
reveals that if new functions~$p$ and~$q$ are defined on~$\Sigma$ by the 
formulae
$$
\aligned
p&=\frac{|b|^2+|a|^2-{1\over9}|a|^4-{9\over{16}}T^2-{27\over{16}}}{1+T^2}\,,\\
q&=\frac{{1\over3}(b{\bar a}^2+{\bar b} a^2)+|a|^2T-{9\over8}T}{1+T^2}\,,\\
\endaligned
\tag15
$$
then equations~(14) imply
$$
\align
3\,dp+{\frac{2q\,dT}{1+T^2}}&=0,\\
3\,dq-{\frac{2p\,dT}{1+T^2}}&=0.\\
\endalign
$$
Thus, the complex valued function~$w$ defined by
$$
w = (p+iq)^3\frac{(1-iT)}{(1+iT)}
$$
must be constant on~$\Sigma$.  Especially, note that the function
$W = p^2+q^2$ is constant on~$\Sigma$, as this remark will be
important in what follows. 

Note that~$w$ can be regarded as the result of composing
the structure function~$(T,a,b):\Sigma\to\bbR\times\bbC\times\bbC$ with 
a $\bbC$@-valued function $\wbold:\bbR\times\bbC\times\bbC\to\bbC$. 
Thus, the constancy of~$w$ can be regarded as the statement that the 
image of the structure function lies in a level surface of~$\wbold$.

\subhead 3.7. A differential ideal \endsubhead
In this section, I will construct a differential ideal on the
manifold~$X=\SL(3,\bbR)\times\bbR\times\bbC\times\bbC$ whose integral
manifolds will correspond to the local and global solutions of the structure 
equations~(14).

Let $\gbold:X\to\SL(3,\bbR)$, $\Tbold:X\to\bbR$, $\abold:X\to\bbC$, 
and $\bbold:X\to\bbC$ denote the projections onto the first, second, third,
and fourth factors of~$X$ respectively.  Set
$$
\spreadmatrixlines{3pt}
\gbold^{-1}\,d\gbold = 
\pmatrix
\phibold^1_1&\phibold^1_2&\phibold^1_3\\
\phibold^2_1&\phibold^2_2&\phibold^2_3\\
\phibold^3_1&\phibold^3_2&\phibold^3_3\\
\endpmatrix.
$$
Then $d\phibold^i_j = -\phibold^i_k\w\phibold^k_j$ and 
$\phibold^1_1+\phibold^2_2+\phibold^3_3=0$, but this is the only linear
relation among these nine $1$@-forms. 

Define $\Ibold$ and $\Jbold$ on~$X$ by~$\Ibold+i\,\Jbold=-(\Tbold+i)\,\abold$, 
set~$\zetabold=\phibold^3_2+i\,\phibold^3_1$, and set
$\rhobold = -\phibold^2_1+\Ibold\,\phibold^3_1+\Jbold\,\phibold^3_2$. 

Now consider the eleven 1@-forms~$\theta_0,\dots,\theta_{10}$ on~$X$ defined 
by the equations
$$
\align
\theta_0 
&=\phibold^1_1-{\ts{1\over3}}(\Ibold\,\phibold^3_2-\Jbold\,\phibold^3_1)\\
\theta_1 
&=\phibold^2_2-{\ts{1\over3}}(\Ibold\,\phibold^3_2-\Jbold\,\phibold^3_1)\\
\theta_2 
&=\phibold^3_3+{\ts{2\over3}}(\Ibold\,\phibold^3_2-\Jbold\,\phibold^3_1)\\
\theta_3 
&=\phibold^1_2+\phibold^2_1\\
\theta_4
&= \phibold^1_3 + \phibold^3_1 - \Tbold\,\phibold^3_2\\
\theta_5
&= \phibold^2_3 + \phibold^3_2 + \Tbold\,\phibold^3_1\\
\theta_6
&=d\Tbold-(\Tbold^2+1)\bigl(\abold\,\zetabold+\bar\abold\,\bar\zetabold\bigr)\\
\theta_7+i\,\theta_8
&=d\abold-i\,\abold\,\rhobold-(\bbold +\abold^2\,\Tbold)\,\zetabold
  +{\ts{3\over2}}\Tbold\,\bar\zetabold\\
\theta_9+i\,\theta_{10}
&=d\bbold-2i\,\bbold\,\rhobold
      -(2\abold\bbold\Tbold -{\ts{2\over9}}\abold^3)\,\zetabold
        +\abold\,\bar\zetabold\,.\\
\endalign
$$
These forms satisfy the relation~$\theta_0+\theta_1+\theta_2 = 0$ but
are otherwise linearly independent.  Thus, they generate a Pfaffian system
~$\Theta$ of rank~10 on~$X$. 

The interest in this Pfaffian system is explained by the following result.

\proclaim{Theorem 6}
Let~$\Sigma$ be a compact, simply-connected $3$@-manifold and let~$\omega$ 
be a generalized Finsler structure on~$\Sigma$ that is projectively flat 
and satisfies~$K=1$.  Let~$g:\Sigma\to\SL(3,\bbR)$ be its canonical immersion,
and let~$(T,a,b):\Sigma\to\bbR\times\bbC\times\bbC$ be its structure function.

Then the map~$(g,T,a,b):\Sigma\to X$ immerses~$\Sigma$ as an integral manifold
of~$\Theta$ and is a double cover onto its image, which is diffeomorphic to~
$\SO(3)$ and which defines a Finsler structure on~$\bbS$. 

Conversely, if~$N\subset X$ is any $3$@-dimensional integral manifold 
of~$\Theta$, then the triple~$\omega=(\phibold^2_1,\phibold^3_1,\phibold^3_2)$ 
defines a generalized Finsler structure on~$N$ that is projectively 
flat and satisfies~$K=1$.  Moreover, if $N$ is compact, then $N$ 
is diffeomorphic to~$\SO(3)$ and the pair~$\bigl(N,\omega\bigr)$ is 
diffeomorphic to a Finsler structure on~$\bbS$ endowed with its canonical 
coframing. 
\endproclaim

\demo{Proof}
Essentially all the pieces of the proof have been assembled in the 
previous sections.  That the map~$(g,T,a,b)$ immerses~$\Sigma$ as
an integral manifold of~$\Theta$ is a consequence of the way that
$g$ was defined to satisfy~$g^{-1}\,dg = \phi$ in \S3.2 and 
the vanishing theorem of \S3.6, which established the equations~(14).
The combination of these results shows that~$(g,T,a,b)^*(\theta_j)
=0$ for $0\le j\le 10$.  Moreover, by the theorem of~\S3.2, it follows
that the map~$g$ double covers onto its image, which defines a
 Finsler structure on~$\bbS$.  Since the structure function~$(T,a,b)$
is defined in terms of the coframing~$\omega$, it must be invariant
under the involution that defines the non-trivial deck transformation of
this cover, so the full map~$(g,T,a,b)$ must also be a double cover
onto its image.

Conversely, since the dimension of~$X$ is 13 and the rank of the Pfaffian
system~$\Theta$ is 10, any integral manifold~$\Sigma\subset X$ of~$\Theta$ 
cannot have dimension more than 3.  Moreover, since the three $1$@-forms
$\phibold^2_1$, $\phibold^3_1$, and $\phibold^3_2$ are linearly independent
modulo~$\Theta$, it follows that, on any $3$@-dimensional integral manifold
$N$ of~$\Theta$, these three $1$@-forms must define a coframing.  The 
vanishing of the forms~$\theta_0,\dots,\theta_5$ on such an $N$
imply that~$\phibold$ has the form of equation~(5) on~$N$. 
Of course this implies that if one defines $\omega_1=\phibold^2_1$, 
$\omega_2=\phibold^3_1$, and $\omega_3=\phibold^3_2$, then~$\omega 
= (\omega_1,\omega_2,\omega_3)$ will satisfy the structure
equations of a generalized Finsler structure on~$N$ that is projectively
flat and satisfies~$K=1$. 

 Finally, since the projection~$g:N\to SL(3,\bbR)$ immerses~$N$ in such
a way that ~$g^{-1}\,dg$ on~$N$ has the form~(5), the arguments of \S3.2
show that, if it is compact, it must be derived from a Finsler structure 
on~$\bbS$ and hence must have the stated topology and geometry. \qed
\enddemo

\subhead 3.8. The uniqueness theorem \endsubhead
The materials are now assembled to prove the uniqueness theorem of
\S3.6.  The crucial observation is that, because $\Theta$ is a Pfaffian
system of rank~10 on the 13@-manifold~$X$, any two integral manifolds
of dimension~3 that intersect must be equal in a neighborhood of any point
of intersection.

\demo{Proof of Theorem~$4$} 
Suppose that $\bigl(\Sigma,\omega\bigr)$ and $\bigl(\Sigma^*,\omega^*\bigr)$
satisfy the hypotheses of the theorem and that $\ubold\in\Sigma$ and
$\ubold^*\in\Sigma^*$ satisfy~$(T,a,b)(\ubold)= (T^*,a^*,b^*)(\ubold^*)$.
Since the canonical immersions~$g:\Sigma\to\SL(3,\bbR)$ and $g^*:\Sigma^*
\to\SL(3,\bbR)$ are only well-defined up to left translation by constants,
they can be uniquely specified by requiring that $g(\ubold) = g^*(\ubold)
= I_3$. 

The theorem of~\S3.7 now implies that each of the maps 
$(g,T,a,b):\Sigma\to X$ and $(g^*,T^*,a^*,b^*):\Sigma^*\to X$ is a double
cover onto a compact 3@-dimensional integral manifold of~$\Theta$.  The 
hypothesis that the structure maps at~$\ubold$ and $\ubold^*$ assume the 
same value implies that these two image integral manifolds intersect in at 
least one point.  Thus, they must be equal, say, 
$(g,T,a,b)(\Sigma) = (g^*,T^*,a^*,b^*)(\Sigma^*) = N$.  It follows that
there is a unique diffeomorphism~$f:\Sigma^*\to\Sigma$ so that 
$f(\ubold^*)=\ubold$ and $(g^*,T^*,a^*,b^*)=(g,T,a,b)\circ f$.  Since,
by construction, the generalized Finsler structures~$\omega$ and $\omega^*$
are pulled back from the canonical one induced on the integral manifold~$N$,
it follows that $f^*(\omega)=\omega^*$, as desired. \qed
\enddemo

\subhead 3.9. The existence theorem \endsubhead
It remains to determine `how many' compact integral manifolds
of~$\Theta$ there are. 

There are some:  The
codimension~$5$ submanifold~$Z\subset X$ defined by the equations
$\Tbold=\abold=\bbold=0$ is an integral manifold of
the five $1$@-forms~$\theta_6,\dots,\theta_{10}$.  On~$Z$, the remaining
$1$-forms in~$\Theta$ reduce to the system~
$\{\phibold^1_1,\,\phibold^2_2,\,\phibold^3_3,\,
\phibold^1_2{+}\phibold^2_1,\,\phibold^3_2{+}\phibold^2_3,\,
\phibold^3_1{+}\phibold^1_3\}$.  The integral manifolds of this latter system 
are the left cosets of~$\SO(3)\subset\SL(3,\bbR)$. 
By construction, these integral manifolds
correspond to Finsler structures on~$\bbS$ that satisfy~$I=J=0$,
i.e., Riemannian metrics of Gauss curvature~1 on~$\bbS$ whose 
geodesics are the standard straight lines. 

The Riemannian solutions may be regarded as the trivial ones.  I will now
show that there are non-trivial ones that are, in a sense, close
to Riemannian.  The first step is the following proposition, whose proof
is a straightforward calculation from the structure equations, so I omit
it.

\proclaim{Proposition 3}
The Pfaffian system~$\Theta$ is Frobenius.  In particular, $X$ has a
foliation~$\cF$ of codimension~$10$ whose leaves are maximal integral 
manifolds of~$\Theta$. \qed
\endproclaim

\demo{Proof of Theorem 5}
Consider the functions~$\pbold$ and $\qbold$ defined on~$\bbR\times\bbC\times
\bbC$ by the formulae
$$
\align
\pbold &= 
\frac{|\bbold|^2+|\abold|^2-{1\over9}|\abold|^4
-{9\over{16}}\Tbold^2-{27\over{16}}}{1+\Tbold^2}\,,\\
\qbold &= 
\frac{{1\over3}(\bbold{\bar \abold}^2+{\bar \bbold} \abold^2)
+|\abold|^2 \Tbold-{9\over8}\Tbold}{1+\Tbold^2}\,.\\
\endalign
$$
A calculation shows that the $1$@-forms~$\theta_{11}$ and $\theta_{12}$
defined by the formulae
$$
\theta_{11} = 3\,d\pbold+{\frac{2\qbold\,d\Tbold}{1+\Tbold^2}}
\qquad\qquad\text{and}\qquad\qquad
\theta_{12} = 3\,d\qbold-{\frac{2\pbold\,d\Tbold}{1+\Tbold^2}}
$$
are linear combinations of the forms~$\theta_6,\dots,\theta_{10}$.  It
follows that the complex valued function~$\wbold$ defined by
$$
\wbold = (\pbold+i\qbold)^3\frac{(1-i\Tbold)}{(1+i\Tbold)}
$$
must be constant on each leaf of~$\cF$. 

In particular, the function $\Wbold = \pbold^2+\qbold^2$ is 
constant on each leaf of~$\cF$.  Now the Taylor series expansion of~
$\Wbold$ in a neighborhood of~$(\Tbold,\abold,\bbold)=(0,0,0)={\bold0}$ 
is of the form
$$
\Wbold = \frac{27}{8}
\left(\frac{27}{8}-\frac{3}{4}\Tbold^2-|\abold|^2-|\bbold|^2\right)
+\dots
$$
where the omitted terms vanish to order at least~3 at~${\bold0}$. 
It follows that ${\bold0}$ is a non-degenerate local maximum of~$\Wbold$ 
on~$\bbR\times\bbC\times\bbC$. In particular, there exists a compact domain~
$D\subset\bbR\times\bbC\times\bbC$ containing ${\bold0}$ in its interior and 
whose smooth boundary is a compact component of a level set of~$\Wbold$.
Because~$\Wbold$ is constant on the leaves of~$\cF$, it follows that
any leaf of~$\cF$ that intersects~$\SL(3,\bbR)\times D\subset X$ must
lie entirely inside~$\SL(3,\bbR)\times D$.

Consider the natural projection~$h:\SL(3,\bbR)\times D\to\SL(3,\bbR)/P$.
(Recall that $P$ is the subgroup of upper triangular matrices with
positive entries on the diagonal.)  The leaves of~$\cF$ are transverse 
to the fibers of~$h$, so $h$ restricts to each leaf~$L$ of~$\cF$ to be
a local diffeomorphism.  In fact, the manifest invariance of~$\Theta$ 
under the natural left action of~$\SL(3,\bbR)$ on~$X$ combined with
the compactness of~$D$ shows that, for any leaf $L$ of~$\cF$ that
lies in~$\SL(3,\bbR)\times D$, every point of~$\SL(3,\bbR)/P\simeq
\SO(3)$ is evenly covered by the local diffeomorphism~$h:L\to\SL(3,\bbR)/P$.
Thus, each such~$L$ must be a covering space of~$\SO(3)$ and hence must
be compact.  By the previous theorem, such a compact leaf must be
diffeomorphic to~$\SO(3)$ and hence the map~$h$ restricts to each leaf~$L$
to be a diffeomorphism onto~$\SL(3,\bbR)/P$.

 For any $(T_0,a_0,b_0)\in D$, a generalized Finsler structure~$\omega$
that is projectively flat and satisfies~$K=1$ and whose structure function
assumes the value~$(T_0,a_0,b_0)$ can now be constructed on~$S^3$ by simply
taking the universal cover of the leaf~$L$ of~$\cF$ that passes through
$(I_3,T_0,a_0,b_0)\in\SL(3,\bbR)\times D$. \qed
\enddemo

\remark{Remark}  Not all of the leaves of~$\cF$ on~$X$ are compact. 
The topology of the non-compact leaves is yet to be determined,
but for some information, see~\S4.
\endremark

\head 4. Connections with the treatment of Funk \endhead

In this final section, I want to explain the relationship of the 
preceding sections with the earlier results of Funk~\cite{Fu2,Fu3}
on Finsler metrics on the plane whose geodesics are the straight lines
and whose curvature satisfies~$K=1$.

In fact, Funk found a complete local classification, though this was
not apparent to me at first, as I found his results hard to understand 
and his arguments hard to follow. In this section, I will give a discussion 
of his results in language that will ease comparison with the other sections 
of this paper and then go on to use these results to write down explicit 
global examples.

\subhead 4.1. Funk's results \endsubhead
One of the difficulties of reading Funk's work for global implications
is that he works on the affine plane~$\bbR^2$ rather than on the natural
global object~$\bbR\bbP^2$ or its double cover, which I identify as $\bbS$
as in~\S2.2.  Since I will be working on~$\bbS$ (both locally
and globally), his results require some translation before they can be
compared with those of the previous sections.

\subsubhead 4.1.1. Projectively parametrized lines in~$\bbS$ \endsubsubhead
Recall the notation of~\S2.2, where~$V$ is identified with~$\bbR^3$
and~$\bbS$ (diffeomorphic to~$S^2$) is the space of oriented lines
through the origin in~$V$.  I will fix the standard volume form on~$\bbR^3$,
i.e., the standard identification of ~$\Lambda^3(V)$ with $\bbR$.  Thus,
for any three vectors~$v_0$, $v_1$, and $v_2$ in~$V$, the wedge product
$v_0\w v_1\w v_2$ will be treated as a number and this identification is
invariant under the natural action of~$\SL(3,\bbR)$ on~$V$.

Given an oriented $2$@-dimensional subspace~
$E\subset V$, let~$\vbold = (v_0,v_1)$ be an oriented basis of~$E$.  Then the 
oriented line~$[E]$ in~$\bbS$ is defined as the oriented curve parametrized
by the map $\gamma_\vbold:S^1\to\bbS$ defined by the formula
$$
\gamma_\vbold(s) = \bigl[\>\cos s\>\>v_0 + \sin s\>\>v_1\>\bigr]
$$
together with the convention that~$S^1$ be oriented so that~$ds$ is a 
positive $1$@-form on~$S^1$. The choice of oriented basis~$\vbold$ of~$E$ 
affects this parametrization but any two such choices will yield the same 
image line with the same orientation. 

A choice of oriented basis~$\vbold$ also defines a linear 
functional~$\alpha_\vbold:V\to\bbR$ by the rule~$\alpha_\vbold(v) 
= v_0\w v_1\w v$.  The resulting oriented line~$[\alpha_\vbold]\in\bbS^*$ 
depends only on the oriented plane~$E$, not on the choice of oriented basis~
$\vbold$, so I will simply identify the points of~$\bbS^*$ with 
the space of oriented lines in~$\bbS$ via the identification~$[E]
=[\alpha_\vbold]$.

\subsubhead 4.1.2. Projectively parametrized geodesics \endsubsubhead
An open domain~$\bbD\subset\bbS$ will be said to be {\it convex\/} if its
intersection with each line in~$\bbS$ is connected.  If ~$\bbD$ is
convex, then I will let~$\bbD^*\subset\bbS^*$ denote the set of oriented
lines in~$\bbS$ whose intersection with~$\bbD$ is non-empty.  (Note that
it is {\it not\/} generally true that ~$\bbD_1^*=\bbD_2^*$ implies
that $\bbD_1 = \bbD_2$ unless $\bbD_1^*\not=\bbS^*$.)

A Finsler structure~$\Sigma_\bbD\subset T\bbD$ will be said to have
{\it linear geodesics\/} if each of its (oriented) geodesics is of
the form~$[E]\,\smcap\,\bbD$ for some (unique) oriented 2@-plane~$E$. 

I can now state one of Funk's results:

\proclaim{Theorem 7 (Funk)}
Let~$\bbD\subset\bbS$ be a convex domain in~$\bbS$ and suppose that there is
a Finsler structure~$\Sigma_\bbD$ on~$\bbD$ with linear geodesics and whose 
curvature satisfies $K=1$.  Then, for every oriented line~$[E]\in\bbD^*$, 
there exists an oriented basis~$\vbold=(v_0,v_1)$ of~$[E]$ so that the 
parametrization~$\gamma_\vbold$ has unit speed \rom(i.e., is a 
$\Sigma_\bbD$@-curve\rom).
\endproclaim

\demo{Remark}
 Funk's proof is based on the characterization of the Finsler-Gauss curvature 
as the $0$@-th order term in the self-adjoint form of the Jacobi equation 
for variation of geodesics.  It is straightforward, so  I will not reproduce 
it here, see~\cite{Fu2}. \qed
\enddemo

\subsubhead 4.1.3. The space of oriented metric lines \endsubsubhead
Theorem~7 suggests looking at the problem of determining when two
oriented bases of a 2-plane~$E$ induce the same metric and orientation
on the line~$[E]$, and this is what I will do next. 

Suppose that $\vbold=(v_0,v_1)$ and $\wbold=(w_0,w_1)$ be two oriented bases 
of~$E$ with the property that the two parametrizations~$\gamma_\vbold$ and~ 
$\gamma_\wbold$ induce the same metric (and orientation) on the line~$[E]$. 
Then there is a constant~$s_0$ so that ~$\gamma_\vbold(s) = 
\gamma_\wbold(s{-}s_0)$.  I.e., for all~$s$,
$$
\bigl[\>\cos s\>v_0 + \sin s\>v_1\>\bigr]
=\bigl[\>\cos(s{-}s_0)\>w_0 + \sin(s{-}s_0)\>w_1\>\bigr],
$$
which can only hold if there exists a positive real number~$r$
so that, for all~$s$,
$$
\bigl(\>\cos s\>v_0 + \sin s\>v_1\>\bigr)
=r\,\bigl(\>\cos(s{-}s_0)\>w_0 + \sin(s{-}s_0)\>w_1\>\bigr).
$$
This, in turn, is equivalent to
$$
v_0 + i\,v_1 = re^{is_0}\,(w_0+i\,w_1)\,.
$$
In particular, the points~$\lb\,v_0{+}i\,v_1\,\rb$ and 
$\lb\,w_0{+}i\,w_1\,\rb$ in~$\bbC\bbP^2 = \bbP(V\otimes\bbC)$ are equal. 

Now, let~$\bbR\bbP^2\subset\bbC\bbP^2$ denote the set of real points.  Given
any point~$z$ in~$\bbC\bbP^2\setminus\bbR\bbP^2$, it can be represented
in the form~$z=\lb\,v_0{+}i\,v_1\,\rb$ for some linearly independent (real) 
vectors~$v_0$ and~$v_1$ in~$V$.  The plane~$E_z$ spanned by the pair
$\vbold=(v_0,v_1)$ and the orientation for which this is an oriented basis
are independent of the specific choice of~$v_0$ and $v_1$ satisfying
$z=\lb\,v_0{+}i\,v_1\,\rb$.  Likewise, 
the metric on~$[E_z]$ for which~$\gamma_\vbold$ is a unit speed 
parametrization does not depend on this choice, but only on~$z$. 

In this way, the open 4-manifold~$\bbC\bbP^2\setminus\bbR\bbP^2$ can be
regarded as a space of oriented Riemannian metrics on lines in~$\bbS$.
The map~$\tau:\bbC\bbP^2\setminus\bbR\bbP^2\to\bbS^*$ defined by
$\tau(z) = [\alpha_\vbold]$ is a smooth submersion whose fiber over~$[E]
\in\bbS^*$ consists of a two-parameter family of Riemannian metrics on
the oriented line~$[E]$.

Theorem~7 can now be reformulated:

\proclaim{Theorem 7$'$}
Let~$\bbD\subset\bbS$ be a convex domain in~$\bbS$ and suppose that there is
a Finsler structure~$\Sigma_\bbD$ on~$\bbD$ with linear geodesics and whose 
curvature satisfies~$K=1$.  Then there exists a unique section~$\sigma:\bbD^*
\to\bbC\bbP^2\setminus\bbR\bbP^2$ of the bundle 
$\tau:\bbC\bbP^2\setminus\bbR\bbP^2\to\bbS^*$ so that, for each~$[E]\in\bbD^*$,
the metric~$\sigma\bigl([E]\bigr)$ on~$[E]$ agrees with the metric 
induced on~$[E]\,\smcap\,\bbD$ by~$\Sigma_\bbD$. \qed
\endproclaim

\remark{Remark}
The section~$\sigma$ will be called the {\it canonical section\/} associated 
to~$\Sigma_\bbD$. 
\endremark

\subsubhead 4.1.4. Local Finsler structures and complex curves\endsubsubhead
The advantage of recasting Funk's result in this language is that it simplifies 
the statement of Funk's local characterization of these Finsler structures:

\proclaim{Theorem 8}
Let~$\bbD\subset\bbS$ be a convex domain in~$\bbS$ and suppose that there is
a Finsler structure~$\Sigma_\bbD$ on~$\bbD$ with linear geodesics and whose 
curvature satisfies $K=1$.  Then the image of the canonical section~
$\sigma:\bbD^*\to\bbC\bbP^2\setminus\bbR\bbP^2$ is a complex curve in~
$\bbC\bbP^2$.

Conversely, if~$C\subset\bbC\bbP^2{\setminus}\bbR\bbP^2$ is a complex
curve with the property that the map~$\tau:C\to\bbS^*$ is a diffeomorphism
onto its image, then there exists a generalized Finsler structure
~$\Sigma_C\subset T\bbS$ with the property that $\tau(z)$ endowed with
the metric~$z$ is a ~$\Sigma_C$@-curve.  Furthermore, this generalized
 Finsler structure has the lines in~$\bbS$ as geodesics and satisfies
$K=1$. \qed
\endproclaim

\remark{Remarks}
Theorem~8 is only partly equivalent to the corresponding results
of~\cite{Fu2}.  Indeed, it took me quite some time to make the translation
of his results into this form.  What Funk does have is a formula for the
local Finsler structures with linear geodesics and curvature equal to~1
in terms of an arbitrary holomorphic function of one variable, which is
locally the same as a complex curve in $\bbC^2$. 

This should be compared with the constructions of \S\S3.4--5, where 
manipulation of the structure equations of a projectively flat generalized 
 Finsler structure with~$K=1$ leads to the introduction of a complex structure 
on the space of geodesics and then the discovery of a holomorphic cubic 
differential, in terms of which the local geometry can be recovered.

The main advantage of this version of Funk's local characterization theorem 
is that it is fully invariant under the action of the projective group~
$\SL(3,\bbR)$, which Funk's description was not.  Even in the later paper~
\cite{Fu3}, where Funk was concerned with a `global' characterization of
the Riemannian metrics on the plane having constant curvature and linear
geodesics among the Finsler structures on the plane with these properties,
he is seriously hampered by not having a formulation that is invariant
under the full projective group. 
\endremark

\subsubhead 4.1.5. A global classification \endsubsubhead
Combining these results of Funk leads immediately to a characterization of 
the Finsler structures defined on the entire sphere~$\bbS$ with linear 
geodesics and with $K=1$:

\proclaim{Theorem 9}
Let~$\Sigma$ be a Finsler structure on~$\bbS$ with linear geodesics
that satisfies~$K=1$.  The image of the associated canonical section~
$\sigma:\bbS^*\to \bbC\bbP^2{\setminus}\bbR\bbP^2$ is then a smooth conic~$C$ 
\rom(i.e., smooth algebraic curve of degree~$2$\rom) without real points.

Conversely, if~$C\subset\bbC\bbP^2$ is a smooth conic without real points, 
then it is the image of the canonical section of a unique Finsler structure~
$\Sigma_C$ on~$\bbS$ with linear geodesics that satisfies~$K=1$.
\endproclaim

\demo{Proof}
By Theorem~8, the image~$C$ is a complex curve in~$\bbC\bbP^2{\setminus}
\bbR\bbP^2$. Since it is the image of a section of a bundle whose base is 
diffeomorphic to a $2$-sphere, it follows that $C$ must be diffeomorphic to a 
$2$@-sphere, i.e., it must be a smoothly embedded rational curve in~
$\bbC\bbP^2$ without real points.  Now, it is a standard result in algebraic 
geometry~\cite{GH, Chapter~1} that any smoothly embedded rational curve 
in~$\bbC\bbP^2$ is either a line or a smooth conic.  Since any line in~
$\bbC\bbP^2$ must have a real point, it follows that~$C$ must be a smooth 
conic.

Conversely, suppose that~$C\subset\bbC\bbP^2$ is
a smooth conic without real points.  I am going to show that $C$ is the
image of a smooth section of~$\tau$.  To do this, it suffices to
show that~$\tau$ restricts to~$C$ to be a diffeomorphism onto~$\bbS^*$.

Let~$[E]\in\bbS^*$ be an oriented line spanned by 
the oriented $2$-plane~$E\subset V$. Consider the real line~
$\bbP\bigl(E\otimes\bbC\bigr)\subset\bbC\bbP^2$.  Since~$C$ is smooth,
it cannot have this line as a component and so, by Bezout's theorem~
\cite{GH, Chapter 1}, it must either intersect~$C$ transversely in two 
distinct points or else in a single point of tangency.  However, it
turns out that a conic without real points cannot have any real tangent lines, 
so $C$ must intersect each of these real lines~
$\bbP\bigl(E\otimes\bbC\bigr)\subset\bbC\bbP^2$ transversely in two points. 

Now, the fibers of~$\tau$ are precisely the connected components of the real 
lines minus their real points.  It follows that~$\tau$ restricts to~$C$ to 
be a submersion from~$C$ to~$\bbS^*$ that is at most 2-to-1.  Since both~$C$ 
and $\bbS^*$ are compact, $\tau$ must restrict to~$C$ to be a covering space. 
 Finally, since both~$C$ and $\bbS^*$ are diffeomorphic to the 2@-sphere, 
$\tau$ must restrict to~$C$ be a diffeomorphism onto~$\bbS^*$. \qed
\enddemo

\remark{Remark}
The global characterization of smooth rational curves in~$\bbC\bbP^2$ as
plane conics should be thought of as corresponding to the vanishing theorem
of~\S3.6.
\endremark

\subhead 4.2. Explicit formulae \endsubhead
In this section, the space of smooth conics in~$\bbC\bbP^2$ 
without real points will be used to study the geometry
of projectively flat Finsler structures on~$\bbS$ with~$K=1$.  Moreover,
an explicit formula will be given for the Finsler norm on~$\bbS$ induced
by such a conic.

\subsubhead 4.2.1.  Conics without real points \endsubsubhead
The natural action of $\SL(3,\bbR)$ on~$V$ extends to the complexification
of~$V$ and to its projectivization~$\bbC\bbP^2$ as well.  This action
preserves the projectivization of~$V$ itself, namely~$\bbR\bbP^2\subset
\bbC\bbP^2$.  Thus, ~$\SL(3,\bbR)$ acts on the set of conics without real
points.  Say that two such conics are {\it $\bbR$@-equivalent\/} if they are
equivalent under this action of~$\SL(3,\bbR)$.

\proclaim{Proposition 4}  Fix a basis~$x,y,z$ of~$V^*$.
Let $p$~and~$q$ be real numbers satisfying~$|q|\le p<\pi/2$.  Then
the conic~$C_{p,q}\subset\bbC\bbP^2$ defined by the equation
$$
e^{ip}\,x^2 + e^{iq}\,y^2 + e^{-ip}\,z^2 = 0
\tag1
$$
is smooth and without real points.  Any conic~$C\subset\bbC\bbP^2$ 
without real points is $\bbR$@-equivalent to~$C_{p,q}$ for some unique~
$(p,q)\in\bbR^2$ satisfying~$|q|\le p<\pi/2$.
\endproclaim

\demo{Proof}
Let $x$, $y$, $z$, and $(p,q)$ be as in the statement of the proposition.
Write
$$
Q = e^{ip}\,x^2 + e^{iq}\,y^2 + e^{-ip}\,z^2 = Q_1 + i\,Q_2
$$
where~$Q_1$ and~$Q_2$ are real quadratic forms on~$V$.  Then $Q_1$ is
positive definite because of the assumptions on~$p$ and $q$, so $Q$ cannot
vanish on any non-zero element of~$V$.  Thus, the conic~$C_{p,q}$ that
it defines via equation~(1) has no real points.  Since~$Q$ is non-degenerate,
$C_{p,q}$ is smooth.

 For the remainder of the proposition, let~$C\subset\bbC\bbP^2$ be a
conic without real points, defined as the null directions of a complex-valued 
quadratic form~$Q = Q_1 + i\,Q_2$ where $Q_1$ and~$Q_2$ are real quadratic 
forms on~$V$. Note that~$Q$ is uniquely determined by~$C$ up to a complex 
scalar multiple.

A slightly messy but straightforward argument shows that unless~$Q_1$~and~
$Q_2$ are simultaneously diagonalizable they will have a common (real)
null vector, which is impossible since~$C$ has no real points.  Thus, let
$x$,~$y$,~and~$z$ be basis elements of~$V^*$ so that
$$
Q = a_1\,x^2 + a_2\,y^2 + a_3\,z^2
$$
where the~$a_i$ are complex numbers.  The positive real linear combinations
of the~$a_i$ cannot contain $0\in\bbC$ since, otherwise~$Q$ would have
a non-zero real null vector.  By multiplying~$Q$ by a non-zero complex 
number, it can be arranged that the positive linear combinations of the~
$a_i$ consist of the positive real linear combinations of~$\{e^{ip},e^{-ip}\}$
for some unique real number~$p$ satisfying $0\le p<\pi/2$.  By scaling
and rearranging the elements of the basis~$x,y,z$, it can be arranged that
$a_1=e^{ip}$ and $a_3=e^{-ip}$ while~$a_2 = e^{iq}$ for some unique~$q$
satisfying~$-p\le q\le p$.  \qed
\enddemo

\remark{Remarks}
When $|q|<p$, the stabilizer of~$C_{p,q}$ in~$\SL(3,\bbR)$ is isomorphic 
to~$\bbZ_2{\times}\bbZ_2$ and consists of the unimodular transformations
of the form~$(x,y,z)\mapsto(\pm x,\pm y,\pm z)$.  When $0<|q|=p$,
the stabilizer is~$\text{O}(2)$.  When $p=q=0$, the stabilizer is~$\SO(3)$.

The complex-valued quadratic form
$$
Q = e^{ip}\,x^2 + e^{iq}\,y^2 + e^{-ip}\,z^2
$$
constructed in the course of the proof of Proposition~4 is uniquely determined 
by~$C$ up to a positive real multiple.  A quadratic~$Q$ of this form will
be said to be {\it normalized\/}.  One could further require that
$x\w y\w z$ be a unit volume form on~$V$, which would make~$Q$ unique,
but that will not be needed below.
\endremark

\subsubhead 4.2.2.  Finsler structures on~$\bbS$ \endsubsubhead
In most of this article, a Finsler structure has been defined to be a 
hypersurface in the tangent bundle of a surface.  However, a Finsler 
structure~$\Sigma$ can also be specified by its corresponding `norm' on the 
tangent bundle, i.e., the function~$|\cdot|_\Sigma$ on the tangent bundle 
that satisfies~$|\lambda v|_\Sigma=\lambda\,|v|_\Sigma$ for all tangent 
vectors~$v$ and $\lambda\ge0$ and ~$|v|_\Sigma=1$ for all~$v\in\Sigma$. 
(The quotes around `norm' are intended to remind the reader that~$|-v|$
is not generally the same as~$|v|$.) 

In this final section, it will be convenient to describe the Finsler 
structures of Theorem~9 in terms of their corresponding norms.

As discussed in \S4.1, for any pair of vectors~$\vbold
=(v, w)$ in $V = \bbR^3$ with $v\not=0$, a curve~$\gamma_\vbold$ 
can be defined in~$\bbS$ by the formula
$$
\gamma_\vbold(t) = [\,v + t\,w\,].
\tag2
$$
The velocity at~$t=0$ of this curve will be denoted~$[v,\,w]$ and, is, of
course, an element of the vector space~$T_{[v]}\bbS$. 

The identity~$c\,[v,\,w]=[v,\,c\,w]$ holds for all real numbers~$c$. 
Moreover,
$$
[v,\,w] = [a\,v,\,\,a\,w+b\,v]
\tag3
$$
for all real numbers~$a>0$ and~$b$.  Conversely, if $[v',\,w']=[v,\,w]$,
then $[v',\,w'] = [a\,v,\,\,a\,w+b\,v]$ for some real numbers
$a>0$ and $b$. 
 
It follows that a Finsler norm on~$\bbS$ determines and is determined by a 
function~$F$ on~$\bigl(V\setminus\{\bold0\}\bigr)\times V$ that satisfies
$$
 F\bigl(v,\,w\bigr) = F\bigl(a\,v,\,\,a\,w+b\,v\bigr)
\qquad\text{for all real~$a>0$ and $b$}
\tag4
$$
as well as the homogeneity condition~
$$
 F\bigl(v,\,c\,w\bigr) = c\,F\bigl(v,\,w\bigr)
\qquad\text{for all~$c\ge0$.}
\tag5
$$
 
\proclaim{Theorem 10}
Let~$C\subset\bbC\bbP^2$ be a conic without real points and let~$Q$ be a 
normalized quadratic form on~$V\otimes\bbC$ so that
$C = \{\,\lb v\rb\,\mid\,Q(v)=0\}$.  Let the inner product of two vectors~$v$ 
and $w$ with respect to~$Q$ be denoted $v\cdot w$.  Set
$$
 F_C(v,w) 
= \Re\left[\,
\frac
{\sqrt{\,(w\cdot w)(v\cdot v)-(w\cdot v)^2\,}-i\,(v\cdot w) }
{(v\cdot v)}
\,\right].
\tag6
$$
Then~$F_C$ defines the Finsler norm of the Finsler structure on~$\bbS$ with 
linear geodesics and $K=1$ whose canonical section~$\sigma:\bbS^*\to\bbC\bbP^2$
has its image equal to~$C$.
\endproclaim

\remark{Remark}
Before beginning the proof (which will be a calculation), it probably is
a good idea to explain that the quantity inside the radical in the formula 
for~$F_C$ can never be a negative real number and the branch of the complex
square root being used is the one satisfying~$\sqrt1=1$ and having
the negative real axis as its branch locus.  Formula~(6) then defines
a function on~$\bigl(V\setminus\{\bold0\}\bigr)\times V$ satisfying (4)
and~(5).  It will be seen below that this function is positive and smooth away 
from points of the form~$(v,\,c\,v)$.  (These points represent the zero 
section of~$T\bbS$.)

By the remarks after Proposition~4, the normalized $Q$ is only determined 
by~$C$ up to a positive real multiple.   However, replacing~$Q$ by~$r\,Q$ 
for any positive~$r$ does not affect the resulting function~$F_C$.
\endremark

\demo{Proof}
Let~$V_{(2)}\subset V\times V$ be the set of pairs of linearly independent
vectors in~$V$.  This is a connected, open subset of~$V\times V$.

Let $(v,w)$ be an element of~$V_{(2)}$.  Since~$C$ has no real points and
no real tangents, the complexified line spanned by~$v$ and $w$ must 
intersect~$C$ transversely in two non-real points. 

Let~$p=\lb\,\alpha\,v+\beta\,w\,\rb$ be such a point. 
Since it is not real, neither~$\alpha$ nor~$\beta$ can vanish, nor can
the ratio~$\alpha/\beta$ be real.  By multiplying~$\alpha$ and~$\beta$
by an appropriate scalar, it can be assumed that~$\beta = ib$ for some
non-zero real number~$b$.  Then~$\alpha$ cannot be pure imaginary, so by 
dividing both~$\alpha$ and $\beta$ by the real part of~$\alpha$, it can be 
assumed that~$\alpha = 1{+}ia$ for some real number~$a$.  Thus, such a~
$p\in C$ can be uniquely written in the form
$$
p = \lb\,(1{+}ia)\,v + ib\,w\,\rb
$$
for some real numbers~$a$ and~$b\not=0$. 

Now, there are two such points,
$$
\align
p_1 &= \lb\,(1{+}ia_1)\,v + ib_1\,w\,\rb = \lb\,v + i(b_1\,w+a_1\,v)\,\rb\,,\\
p_2 &= \lb\,(1{+}ia_2)\,v + ib_2\,w\,\rb = \lb\,v + i(b_2\,w+a_2\,v)\,\rb\,.\\
\endalign
$$
Since~$\tau(p_i)=[\,v\w(b_i\,w+a_i\,v)\,]=[\,b_i\,v\w w\,]$ and these two
points must represent the two different oriented lines spanned by~$v$~and~$w$,
it follows that~$b_1$ and $b_2$ are always of different signs.  Thus, it
can be arranged that~$b_1<0<b_2$.  Once this is done,~$b_1$ and~$b_2$ 
can be regarded as smooth (in fact, analytic) functions on~$V_{(2)}$. 
Note that, while~$b_2(v,-w)=-b_1(v,w)$, it will not generally be true 
that~$b_1(v,w)=-b_2(v,w)$ unless the conic~$C$ is invariant under conjugation.

By Theorem~9, the Finsler structure~$\Sigma_C$ is defined on~$\bbS$ so 
that the tangent vector~$[v,\,(b_2\,w+a_2\,v)]= b_2\,[v,\,w]$ is a unit 
vector.  To establish Theorem~10, it suffices to show that the formula 
for~$b_2$ is just
$$
b_2(v,w) = \frac1{F_C(v,w)}\,.
$$
This should be just a calculation, but
it turns out that some argument is needed to show that the formula~$F_C$ 
makes sense, particularly as regards the choice of sign inherent in the
definition of the square root.  It is here that having~$Q$ be
normalized is important.

 First, I will show that the quantity under the radical in~(6) is 
never a negative real number and is, in fact, never zero for~$(v,w)\in
V_{(2)}$.  Since~$C$ is without real points and~$Q$ is normalized, 
Proposition~4 implies that there exist non-negative real numbers~$p$~and~$q$ 
and a basis of~$V^*$, i.e., an identification of~$V$ with~$\bbR^3$, so that, 
for any~$v=(x,y,z)$, the formula
$$
Q(v)  = v\cdot v = e^{ip}\,x^2 + e^{iq}\,y^2 + e^{-ip}\,z^2
$$ 
holds.  If $w=(a,b,c)\in V$ is any other vector, then, of course,
$$
w\cdot w = e^{ip}\,a^2 + e^{iq}\,b^2 + e^{-ip}\,c^2
$$ 
and 
$$
w\cdot v  = e^{ip}\,ax + e^{iq}\,by + e^{-ip}\,cz.
$$
Calculation now yields that $(v\cdot v)(w\cdot w) - (w\cdot v)^2$ is equal to 
$$
e^{i(q-p)}\,(yc-zb)^2 + (xc-za)^2 + e^{i(q+p)}\,(xb-ya)^2.
$$
With a slight abuse of notation, denote this quantity by~$Q(v\w w)$. The
above expression for~$Q(v\w w)$ shows that it lies in the wedge~$W_{p,q}$ 
consisting of the non-negative real linear combinations of~$e^{i(q+p)}$,~$1$, 
and~$e^{i(q-p)}$.  Since, by hypothesis, $|q|<p\le \pi/2$, this wedge does not 
contain any negative real numbers. Thus,~$Q(v\w w)$ cannot be a negative real 
number and cannot equal zero unless $(yc-zb)=(xc-za)=(xb-ya)=0$, i.e., unless 
$v$ and $w$ are linearly dependent. In particular, taking the complex square 
root function to be branched along the negative real axis and to satisfy~
$\sqrt1 = 1$, the expression~$\sqrt{Q(v\w w)}$ defines a smooth (in fact, 
analytic) function on~$V_{(2)}$ with values in the right half plane of~$\bbC$.

Now, for given~$(v,w)\in V_{(2)}$, consider the defining equation for~$p_2$. 
This is
$$
0
=\bigl((1{+}ia_2)\,v+ib_2\,w\bigr)\cdot\bigl((1+ia_2)\,v+ib_2\,w\bigr),
$$
Together with the inequality~$b_2>0$.
Expanding this equation, dividing by $v\cdot v\not=0$ and collecting yields
$$
0 = \left(1+ ia_2 + ib_2\,\frac{(v\cdot w)}{(v\cdot v)}\right)^2 
- {b_2}^2\,\frac{Q(v\w w)}{(v\cdot v)^2}.
$$
After multiplying by $-1$, the right hand side of this equation factors
as
$$
\left(
a_2 + b_2\,\left(\frac{(v\cdot w)}{(v\cdot v)}
+i\,\frac{\sqrt{Q(v\w w)}}{(v\cdot v)}\right) - i
\right)
\left(
a_2 + b_2\,\left(\frac{(v\cdot w)}{(v\cdot v)}
-i\,\frac{\sqrt{Q(v\w w)}}{(v\cdot v)}\right) - i
\right).
$$
Each of these factors is a real analytic function on~$V_{(2)}$, 
so one of them must vanish identically.  Accordingly, write
$$
a_2 + b_2\,\left(\frac{(v\cdot w)}{(v\cdot v)}
+i\varepsilon\,\frac{\sqrt{Q(v\w w)}}{(v\cdot v)}\right) = i\,,
\tag7
$$
where~$\varepsilon=\pm1$.  Since~$a_2$ is real and $b_2$ is real and 
positive, the quantity in the parentheses must have positive
imaginary part.  Note that the real and imaginary parts of (7) 
constitute a pair of linearly independent equations for~$a_2$ and $b_2$.

The same argument applied to~$p_1$ yields
$$
a_1 + b_1\,\left(\frac{(v\cdot w)}{(v\cdot v)}
-i\varepsilon\,\frac{\sqrt{Q(v\w w)}}{(v\cdot v)}\right) = i\,.
\tag8
$$
(The opposite sign must hold since~$b_1\not=b_2$.)   Since~$a_1$ is real
and~$b_1$ is real and negative, the 
quantity in the parentheses in~(8) must have negative imaginary part. 

It now follows from~(7) and~(8) and the signs of the~$b_i$ that the real part 
of 
$$
\varepsilon\,\frac{\sqrt{Q(v\w w)}}{(v\cdot v)}
$$
must be positive on~$V_{(2)}$.  Evaluating this expression at
$(v,w) = \bigl((1,0,0),(0,1,0)\bigr)$ then yields $\varepsilon = +1$. 

Thus, 
$$
a_2 + b_2\,\left(\frac{(v\cdot w)}{(v\cdot v)}
+ i\,\frac{\sqrt{Q(v\w w)}}{(v\cdot v)}\right) 
= 
a_1 + b_1\,\left(\frac{(v\cdot w)}{(v\cdot v)}
- i\,\frac{\sqrt{Q(v\w w)}}{(v\cdot v)}\right) 
= i.
$$
 Finally, taking imaginary parts yields
$$
b_2\,\Re\left[\,
\frac{\sqrt{Q(v\w w)}}{(v\cdot v)}-i\frac{(v\cdot w)}{(v\cdot v)}
\,\right] = 1
$$
as was to be shown.
\qed
\enddemo

\remark{Remarks}
In conclusion, here are a few remarks about the implications of Theorem~10. 

 First, for~$F_C$ defined as in~(6),
$$
 F_C(v,w)-F_C(v,-w) = 
2\Im\left[\,\frac{(v\cdot w)}{(v\cdot v)}\,\right]\,.
$$
It follows that, if the Finsler structure is to be symmetric, then 
$(v\cdot w)/(v\cdot v)$ must be real for all non-zero~$v\in V$ and all~$w
\in V$.  Referring to the formula for a normalized~$Q$, this can only be true 
if $p=q=0$, i.e., if~$Q$ is real-valued.  In this case, formula~(6) 
reduces to the standard formula for the norm of a Riemannian metric and, 
of course, the Finsler structure is Riemannian.  Thus, the only 
symmetric Finsler structure among these Finsler structures is the 
Riemannian one.

Second, the proof shows that $\sqrt{Q(v,w)}/(v\cdot v)$ always has positive
real part when $v\w w\not=0$.  In particular, this implies that 
$Q(v\w w)/(v\cdot v)^2$ is never a negative real number and equals zero
only if $v\w w = 0$.  This fact does not appear to be easy to establish
directly, but, using it, the formula for~$F_C$ can be rewritten in the form
$$
 F_C(v,w) 
= \Re\left[\,
\sqrt{\frac{(w\cdot w)(v\cdot v)-(w\cdot v)^2}{(v\cdot v)^2}}
 -i\,\frac{(v\cdot w)}{(v\cdot v)}
\,\right],
\tag9
$$
a form in which the normalization of~$Q$ is irrelevant.

Third, on any particular tangent space~$T_{[v]}\bbS$, the Finsler
norm has the form
$$
\align
 F_{[v]} &= \Re\left(\,\sqrt{P+i\,Q}\,\right) - L\\
        &= \sqrt{\frac{\sqrt{P^2+Q^2} + P}{2}}-L
\endalign
$$
where ~$P$ and $Q$ are real quadratic forms on~$T_{[v]}\bbS$ while $L$ is 
a real linear form on this space.  It follows without difficulty that the
unit vectors of the Finsler structure in each tangent space form an algebraic 
curve of degree~4, or, in certain degenerate cases, of degree~2. 
In particular, these spaces are not Randers spaces, nor, indeed, any of the 
other special types of Finsler surfaces considered in~\cite{Ma2}.

 Fourth, for a fixed basis~$x,y,z$ of~$V^*$, as $(p,q)$ ranges through
the triangle~$|q|\le p<\pi/2$, the corresponding curves~$C_{p,q}$ give rise
to Finsler norms~$F_{p,q}$ on~$\bbS$ that are inequivalent.  This
two parameter family must correspond to the two parameter family of 
inequivalent compact integral manifolds of the Pfaffian system~$\Theta$ 
constructed in~\S3.8. 

It seems that, by some principle of analytic 
continuation, the non-compact integral manifolds must correspond to the
partial Finsler structures on~$\bbS$ induced by conics~$C\subset\bbC\bbP^2$ 
that have one or more real points.  It would be interesting to see
what the formula~$F_C$ gives for such conics.

 Fifth, the values~$(p,q)$ where $q=\pm p\not=0$ correspond to curves~$C_{p,q}$
that have a one-parameter symmetry group in~$\SL(3,\bbR)$.  The corresponding
 Finsler structures must also be invariant under this subgroup and hence are
 Finsler surfaces of rotation.  Thus, there exists a one parameter family of 
rotationally invariant, projectively flat Finsler structures with~$K=1$.

Even when~$|q|<p$, there is still a~$\bbZ_2{\times}\bbZ_2{\times}\bbZ_2$ 
symmetry group of the Finsler structure~$F_{p,q}$ (with an index~2 subgroup
consisting of the orientation preserving symmetries). 
These correspond to reflectional symmetries in three `axis' geodesics.

 Finally, since $F_C(v,w) = F_C(-v,-w)$ it follows that each of these Finsler
structures is invariant under the antipodal map on~$\bbS$ and hence
descends to a well-defined Finsler structure on~$\bbR\bbP^2$ with linear
geodesics and Finsler-Gauss curvature equal to~$1$.
\endremark

\enddocument